\documentclass[10pt]{article}

\usepackage{amsmath}
\usepackage{amssymb}
\usepackage{indentfirst}
\usepackage{graphics} 
\usepackage{color}

\makeatletter

\@addtoreset{equation}{section} 
\makeatother
\def\<{\langle}
\def\>{\rangle}

\newtheorem{lem}{Lemma}[section]
\newtheorem{theo}{Theorem}[section]
\newtheorem{rem}{Remark}[section]
\newtheorem{pro}{Proposition}[section]

\makeatletter
   
   \@addtoreset{equation}{section}
\makeatother

\begin{document}
\title{\bf  Double diffusion structure of logarithmically damped wave equations with a small parameter}

\author{ Alessandra Piske\thanks {Corresponding author: alessandrapiske@gmail.com}  \, and Ruy Coimbra Char\~ao\thanks{ ruy.charao@ufsc.br} \;
\\{ \small Graduate Program in Pure and Applied Mathematics}  
 \\{\small Department of Mathematics} \\{\small Federal University of Santa Catarina} \\ {\small 88040-270, Florianopolis, Brazil,} 
\\
and\\Ryo Ikehata\thanks{ikehatar@hiroshima-u.ac.jp} \\ {\small Department of Mathematics}\\ {\small Division of Educational Sciences}\\ {\small Graduate School of Humanities and Social Sciences} \\ {\small Hiroshima University} \\ {\small Higashi-Hiroshima 739-8524, Japan}}
\date{}
\maketitle
\begin{abstract}
We consider a wave equation with a nonlocal logarithmic damping depending on a small parameter $0 < \theta < \frac{1}{2}$. This research is a counter part of that was initiated by Char\~ao-D'Abbicco-Ikehata considered in \cite{CDI} for the large parameter case $\theta > \frac{1}{2}$. We study the Cauchy problem for this model in ${\bf R}^{n}$ to the case $\theta \in (0,\frac{1}{2})$, and we obtain an asymptotic profile and optimal estimates in time of solutions as $t \to \infty$ in $L^{2}$-sense. An important discovery in this research is that in the case when $n = 1$, we can present a threshold $\theta^{*} = \frac{1}{4}$ of the parameter $\theta \in (0,\frac{1}{2})$ such that the solution of the Cauchy problem decays with some optimal rate for $\theta \in (0,\theta^{*})$ as $t \to \infty$, while the $L^{2}$-norm of the corresponding solution never decays for $\theta \in [\theta^{*}, \frac{1}{2})$, and in particular, in the case $\theta \in [\theta^{*},\frac{1}{2})$ it shows an infinite time $L^{2}$-blow up of the corresponding solutions. The former (i.e., $\theta \in (0,\theta^{*})$ case) indicates an usual diffusion phenomenon, while the latter (i.e., $\theta \in [\theta^{*},\frac{1}{2})$ case) implies, so to speak, a singular diffusion phenomenon. Such a singular diffusion in the one dimensional case is a quite novel phenomenon discovered through our new model produced by logarithmic damping with a small parameter $\theta$. It might be already prepared in the usual structural damping case such as $(-\Delta)^{\theta}u_{t}$ with $\theta \in (0,1/2)$, however unfortunately nobody has ever just pointed out even in the structural damping case. 
\end{abstract}
\section{Introduction}
\footnote[0]{Keywords and Phrases: Wave equation; logarithmic damping; small parameter; $L^{2}$-decay; asymptotic profile; optimal estimates; double diffusion phenomenon, singularity.}
\footnote[0]{2010 Mathematics Subject Classification. Primary 35B40; 35L05; Secondary 35B20, 35R12, 35S05.}

We consider in this work  the  dissipative wave equation based on an operator $L_{\theta}$, that combines the composition of  logarithm function  with the Laplace operator,  as follows:
\begin{align}
& u_{tt} -\Delta u + L_{\theta}u_t  = 0,\ \ \ (t,x)\in (0,\infty)\times {\bf R}^{n},\label{eqn}\\
& u(0,x)= u_0(x), \quad  u_{t}(0,x)= u_{1}(x),\ \ \ x\in{\bf R}^{n} ,\label{initial}
\end{align}
where the linear operator  
\[L_{\theta}: D(L) \subset L^{2}({\bf R}^{n}) \to L^{2}({\bf R}^{n}), \quad \theta>0,\]
is defined as follows: 
\[D(L_{\theta}) := \left\{f \in L^{2}({\bf R}^{n}) \,\bigm|\,\int_{{\bf R}^{n}}(\log(1+\vert\xi\vert^{2\theta}))^{2}\vert\hat{f}(\xi)\vert^{2}d\xi < +\infty\right\},\]
and for $f \in D(L_{\theta})$,  
\[(L_{\theta}f) (x) := {\cal F}_{\xi\to x}^{-1}\left(\log (1+\vert\xi\vert^{2\theta})\hat{f}(\xi)\right)(x).\]
\noindent
Here, one has just denoted the Fourier transform ${\cal F}_{x\to\xi}(f)(\xi)$ of $f(x)$ by 
\[{\cal F}_{x\to\xi}(f)(\xi) = \hat{f}(\xi) := \displaystyle{\int_{{\bf R}^{n}}}e^{-ix\cdot\xi}f(x)dx, \quad \xi \in {\bf R}^n,\]
as usual with $i := \sqrt{-1}$, and ${\cal F}_{\xi\to x}^{-1}$ expresses its inverse Fourier transform. 
\noindent
Since the operator $L_{\theta}$ is non-negative and self-adjoint in $L^{2}({\bf R}^{n})$ (see \cite{Log-damping}), the square root 
$$L_{\theta}^{1/2}: D(L_{\theta}^{1/2}) \subset L^{2}({\bf R}^{n}) \to L^{2}({\bf R}^{n})$$
can be defined, and is also nonnegative and self-adjoint with its domain
\[D(L_{\theta}^{1/2}) = \left\{f \in L^{2}({\bf R}^{n}) \,\bigm|\,\int_{{\bf R}^{n}}\log(1+\vert\xi\vert^{2\theta })\vert\hat{f}(\xi)\vert^{2}d\xi < +\infty\right\}.\] 
Note that $D(L_{\theta}^{1/2})$ becomes Hilbert space with its graph norm
$$\Vert v\Vert_{D(L_{\theta}^{1/2})} := \left(\Vert v\Vert^{2} + \Vert L_{\theta}^{1/2}v\Vert^{2}\right)^{1/2},$$
where to simplify the notation  we define the  $L^{2}({\bf R}^{n})$-norm  by
\[\Vert\cdot\Vert := \Vert\cdot\Vert_{L^{2}({\bf R}^{n})}.\]

We  also note that  
$$H^{s}({\bf R}^{n}) \hookrightarrow D(L_{\theta}^{1/2}) \hookrightarrow  L^{2}({\bf R}^{n})$$
for any $ s > 0$. Symbolically writing, one can see
\[L_{\theta}= \log(I+(-\Delta)^{\theta}),\]
where $\Delta$ is the usual Laplace operator defined on $H^2({\bf R}^n)$.
\noindent
Since $L_{\theta}$ is constructed by a nonnegative-valued multiplication operator, it is nonnegative and self-adjoint in $L^{2}({\bf R}^{n})$. Then, by a similar argument to \cite[Proposition 2.1]{ITY} based on the Lumer-Phillips Theorem one can find that the problem (1.1)-(1.2) has a unique mild solution
\[u \in C([0,\infty);H^{1}({\bf R}^{n})) \cap C^{1}([0,\infty);L^{2}({\bf R}^{n}))\]
for each $[u_{0},u_{1}] \in H^{1}({\bf R}^{n})\times L^{2}({\bf R}^{n})$, and the associated energy identity holds 
\begin{equation}\label{energy}
E_{u}(t)  +  \int_0^t\Vert L_{\theta}^{1/2}u_t(s,\cdot)\Vert^{2}ds=  E_{u}(0), \quad t>0,
\end{equation} 
where
\[
E_{u}(t) := \frac{1}{2}\left(\Vert u_{t}(t,\cdot)\Vert_{L^{2}}^{2} + \Vert  \nabla u(t,\cdot)\Vert_{L^{2}}^{2}\right).
\]
The identity  \eqref{energy} implies that the total energy is a non increasing function in time because of the existence of the  dissipative term $  L_{\theta}u_{t}$.  \\

We first try to review some known historical results (basically restricting to linear equations). \\

On the strongly damped wave equation such that
\begin{equation}\label{ike-13}
u_{tt} -\Delta u -\Delta u_t  = 0,\ \ \ (t,x)\in (0,\infty)\times {\bf R}^{n},
\end{equation}
one can mention to the celebrated papers \cite{S} and \cite{Po} concerning the $L^{p}$-$L^{q}$ estimates of solutions to the equation \eqref{ike-13}, and all related researches have their origin in there. After \cite{S} and \cite{Po}, a next important topic concerning \eqref{ike-13} is extensively studied by \cite{ITY}, \cite{I-14} and \cite{IO}, that is, they investigate asymptotic profiles, and optimal estimates of the $L^{2}$-norm of solutions as $t \to \infty$. By measuring the solution in terms of $L^{2}$-norm, a singularity near $0$-frequency region of the solution to problem \eqref{ike-13} can be captured precisely. For a higher order asymptotic expansion as $t \to \infty$ of the solution to \eqref{ike-13} it can be derived by \cite{Mi} recently. A higher order asymptotic expansion in time of the square of $L^{2}$-norm of solutions to the equation \eqref{ike-13} can be precisely obtained by \cite{BV-1, BV-2}. In this connection, it should be mentioned that a critical exponent problem to the semi-linear equation:  
\begin{equation*}
u_{tt} -\Delta u -\Delta u_t  = \vert u\vert^{p},\ \ \ (t,x)\in (0,\infty)\times {\bf R}^{n}
\end{equation*}
can be studied by \cite{DR} based on the $L^{p}$-$L^{q}$ estimates due to \cite{S}. Unfortunately, it seems still open to determine the critical exponent $p^{*}$ of the power $p > 1$ of the nonlinearity. As for small data global existence and blow-up results to the strongly damped waves with two different kinds of nonlinearity $\vert u\vert^{p} + \vert u_{t}\vert^{q}$ or (simply) $\vert u_{t}\vert^{q}$, one can cite \cite{CF, CF-2}. The problems \cite{CF, CF-2} are considered in an exterior domain. \\

On the other hand, recently a more general fractionally damped wave equations such that  
\begin{equation}\label{ike-15}
u_{tt} + (-\Delta)^{\sigma}u + (-\Delta)^{\theta}u_t  = 0,\ \ \ (t,x)\in (0,\infty)\times {\bf R}^{n},
\end{equation}
are studied intermittently, however, a research in the case of $\theta = 0$ has initiated to the paper by \cite{K} in 2000, there the author has captured a self-similar profile in asymptotic sense as $t \to \infty$ (in fact, the author treats more general semi-linear problems). A generalized version of \cite{K} has been just published recently in \cite{FIW}. One can also cite \cite{DE-1, DE-2, DE-3}, \cite{DAO}, \cite{DR}, \cite{Y}, \cite{PMR} with $\sigma = 1$ about the study for precise asymptotic profiles and/or critical exponent of the nonlinear problems. $L^{p}$-$L^{q}$ estimates and/or asymptotic profiles in the case of $\frac{1}{2} \leq \theta < 1$ to the linear equation \eqref{ike-15} with $\sigma = 1$ can be considered by \cite{DEP}, \cite{PMR}, \cite{IT} and \cite{NR}. \\

Under these observations on the equation \eqref{ike-15}, quite recently the authors in \cite{Log-damping} have presented a new type of wave equations with a non-local logarithmic damping:      
\begin{equation}\label{ike-16}
u_{tt} -\Delta u +\log(I -\Delta)u_t  = 0,\ \ \ (t,x)\in (0,\infty)\times {\bf R}^{n},
\end{equation}
where $I$ is an identity operator, and the authors investigated asymptotic profiles and optimal decay rates of solutions to problem \eqref{ike-16}. In this connection, one can understand that the study by \cite{Log-damping} has a close relation to some properties of hypergeometric functions. As consequence, one can know that the solution to \eqref{ike-16} has a similar property to that of \eqref{ike-13}. As a next natural problem to \eqref{ike-16}, the authors in \cite{CDI} study the following generalized equation with parameter $\theta > 0$:  
\begin{equation}\label{ike-17}
u_{tt} -\Delta u +\log(I +(-\Delta)^{\theta})u_t  = 0,\ \ \ (t,x)\in (0,\infty)\times {\bf R}^{n}.
\end{equation}
In fact, the authors in \cite{CDI} treat the case of $\frac{1}{2} < \theta < 1$, and investigate asymptotic profiles and optimal decay rates for the solution itself and the total energy. A general theory on hypergeometric functions can be also effectively used in the proof. The profiles include an oscillation property coming from the condition $\theta > \frac{1}{2}$. In this connection, even if $\theta$ is large enough, say $\theta > 1$, the equation \eqref{ike-17} does not have any regularity-loss structures, while in the case when $\theta > 1$ the equation \eqref{ike-15} with $\sigma = 1$ has such a regularity-loss structure as was pointed out in \cite{II, FIM}. A research in \cite{II} has been motivated by an interesting paper by \cite{GGH}. By the way, such a regularity loss structure has been first discovered by S. Kawashima through the research on Timoshenko systems (e.g.,\cite{IHK}). Additionally, a similar study on asymptotic wave-like property of the solution to the equation   
\begin{equation*}
u_{tt} + \log(I -\Delta)u  + \log(I -\Delta)u_t  = 0,\ \ \ (t,x)\in (0,\infty)\times {\bf R}^{n}.
\end{equation*}
can be deeply investigated in \cite{CAI}.\\

Here, we mention an answer to the question: why  do we study the wave equation with a log-damping term in a series of our papers? This comes from an observation below. \\
At first, notice that the properties of the equation such that
\begin{equation}\label{ike-100}
u_{tt}-\Delta u + (-\Delta)^{\theta}u_{t} = 0
\end{equation} 
can be divided into two parts:\\
\noindent
{\rm (a)}\,$0 \leq \theta \leq 1$ $\Rightarrow$ the equation \eqref{ike-100} does not have a regularity-loss structure, and this fact is well-studied.\\
{\rm (b)}\, $\theta > 1$ $\Rightarrow$ the equation \eqref{ike-100} has a regularity-loss structure coming from a priori estimates in the high frequency region (cf., \cite{II}) as is already mentioned.\\
By comparing the equations \eqref{ike-100} with $\theta = 1$ and \eqref{eqn} with $\theta > 0$, one notices
\begin{equation*}
\lim_{\vert\xi\vert \to \infty}\frac{\vert\xi\vert^{2}}{\log(1+\vert\xi\vert^{2\theta})} = \lim_{r \to \infty}\frac{r}{\log(1+r^{\theta})} = +\infty,
\end{equation*} 
for any $\theta > 0$. This implies that the asymptotic behavior in the high frequency region of the solution to the equation \eqref{eqn} is still under control of that of the equation \eqref{ike-100} with $\theta = 1$ (i.e., strongly damped wave case). As a result, the equation \eqref{eqn} does not have any regularity-loss structures for all $\theta > 0$ (even if $\theta$ is big enough). This is one of our merit to study the equation \eqref{eqn} with $\theta > 0$.\\

By the way, it should be noticed that general theories including semigroup approaches concerning the logarithmic Laplacian can be studied in detail by a series of papers and books due to Amann \cite[Chapter III, p. 152]{A}, Nollau \cite{N}, Reed-Simon \cite[p. 317]{RS}, and Weilenmann \cite{WJ}, Chen-Weth \cite{CW}, and in particular, one can refer to an interesting recent published paper \cite{B}, which studies the logarithmic wave operator in terms of the logarithmic negative Dirichlet Laplacian defined in a bounded smooth domain. \\
   
Our main goal in this paper is to find an asymptotic profile of solutions in the $L^{2}$-framework to problem \eqref{eqn}-\eqref{initial} in the case when (ideally speaking) $0 < \theta < \displaystyle{\frac{1}{2}}$, and is to apply them to investigate the optimal  decay rates, depending on the dimension $n$ and the parameter $\theta$, of solutions to problem \eqref{eqn}-\eqref{initial}. The case $0 < \theta < \frac{1}{2}$ is the missing one in the previous researches. Our interest to study the equation \eqref{eqn} is only from a pure mathematical point of view, and without loss of generality we can assume that the initial amplitude $u_0=0$ when one concentrates only on capturing the leading term as time goes to infinity.\\

Now, we introduce the asymptotic profile as $t \to \infty$ of the solutions to problem \eqref{eqn}-\eqref{initial} with $u_{0} = 0$:
\begin{equation*}
\varphi (t,\xi) := \frac{e^{-\frac{|\xi|^2}{\log(1+|\xi|^{2\theta})}t}}{\log (1+|\xi|^{2\theta})} P_1 - \frac{e^{-\log(1+|\xi|^{2\theta})t}}{\log(1+|\xi|^{2\theta})}P_1 =: P_{1}\left(\varphi_{1}(t,\xi) - \varphi_{2}(t,\xi)\right),
\end{equation*}
where the $0$th-moment of the initial velocity $P_{1} \in {\bf R}$ is defined by 
\[P_{1} := \int_{{\bf R}^{n}}u_{1}(x)dx.\]
Then, our main result reads as follows. It should be strongly mentioned that the case of $n = 1$ and $ \theta \geq \frac{1}{4} $ can be included in the result below. This part will be essential in our paper.
\begin{theo}\label{maintheorem}
\begin{itemize}
\item[\rm{(I).}] Let $ n = 1 $,  $ 0 < \theta \leq \displaystyle{\frac{1}{3}}$, and $ u_1 \in L^{1,2\theta}({\bf R})\cap L^{2}({\bf R}) $. Then it holds that 
\begin{align*}
\Vert u(t,\cdot) - {\cal F}_{\xi\rightarrow x}^{-1}(\varphi (t,\xi))(\cdot)\Vert_{L^{2}}  \leq \left\{\begin{matrix}
C(\|u_1\|_{1}+\|u_1\|_{L^{1,2\theta}}) \left ( t^{-\frac{1}{4(1-\theta)}} + \displaystyle{\frac{1}{\sqrt{\theta}}} t^{-\frac{1}{4\theta}}\right ), & \text{ if } 0< \theta \leq \frac{1}{6}, \\ 
 C(\|u_1\|_{1}+\|u_1\|_{L^{1,2\theta}}) \left ( t^{-\frac{1}{4(1-\theta)}} + \displaystyle{\frac{1}{\sqrt{\theta}}} t^{-\frac{\frac{5}{3}-4\theta}{4\theta}} \right ), & \text{ if } \frac{1}{6} < \theta \leq \frac{1}{3}
\end{matrix}\right.  
\end{align*}
for $t\gg 1$, where $u(t,x)$ is a unique solution to problem {\rm (1.1)-(1.2)} with $u_{0} = 0$. 

\item[\rm{(II).}] Let $ n\geq 2 $,  $ 0 < \theta \leq \displaystyle{\frac{5}{12}}$, and $ u_1 \in L^{1,2\theta}({\bf R}^{n})\cap L^{2}({\bf R}^{n}) $. Then  it holds that 
\begin{align*}
\Vert u(t,\cdot) - {\cal F}_{\xi\rightarrow x}^{-1}(\varphi (t,\xi))(\cdot)\Vert_{L^{2}} \leq \left\{\begin{matrix}
C(\|u_1\|_{1}+\|u_1\|_{L^{1,2\theta}}) \left ( t^{-\frac{n}{4(1-\theta)}} + \displaystyle{\frac{1}{\sqrt{\theta}}} t^{-\frac{n}{4\theta}}\right ), & \text{ if } 0< \theta \leq \frac{1}{6}, \\ 
 C(\|u_1\|_{1}+\|u_1\|_{L^{1,2\theta}}) \left ( t^{-\frac{n}{4(1-\theta)}} + \displaystyle{\frac{1}{\sqrt{\theta}}} t^{-\frac{n-4\theta +\frac{2}{3}}{4\theta}} \right ), & \text{ if } \frac{1}{6} < \theta \leq \frac{1}{3},\\ 
C(\|u_1\|_{1}+\|u_1\|_{L^{1,2\theta}}) \left ( t^{-\frac{n-1}{4(1-\theta)}} + \displaystyle{\frac{1}{\sqrt{\theta}}} t^{-\frac{n-1}{4\theta}} \right ),  & \text{ if } \frac{1}{3} < \theta \leq \frac{5}{12}
\end{matrix}\right.
\end{align*}
for  $t \gg 1$, where $u(t,x)$ is a unique solution to problem {\rm (1.1)-(1.2)} with $u_{0} = 0$.  
\end{itemize}

\end{theo}\begin{rem}{\rm 
In the results of Theorem \ref{maintheorem}, one can notice the coefficient $1/\sqrt{\theta}$ in front of each final estimates. By observing this coefficient, one may conclude that we have captured the unique nature for the $\log$-damping (or fractional damping) with parameter $\theta > 0$. This property can be found by searching the leading term more precisely than previous researches .}
\end{rem}
\begin{rem}
{\rm It follows from Theorem \ref{maintheorem} that $\hat{u}(t,\xi) \sim P_{1}\left(\varphi_{1}(t,\xi) - \varphi_{2}(t,\xi)\right)$ in $L^{2}({\bf R}_{\xi}^{n})$ as $t \to \infty$. It is important to notice that $\varphi_{1}(t,\xi)$ and   $\varphi_{2}(t,\xi)$ are exact solutions of the first order in time equations in the Fourier space, respectively:
\[-\Delta v + L_{\theta} v_{t} = 0,\]
and
\[L_{\theta} v + v_{t} = 0.\]
In some sense, the solution to problem (1.1)-(1.2) with small parameters $\theta \in (0,1/2)$ has a double diffusion phenomenon. This kind of important double diffusion phenomenon has been first discovered by D'Abbicco-Ebert \cite{DE-1} to the equation \eqref{ike-100} with $\theta \in (0,1/2)$. Theorem \ref{maintheorem} corresponds to that of \cite[Theorem 2]{DE-3}. We find that (1.1)-(1.2) has a similar property to it. While, in the case when $n \geq 2$ and $\theta \in (0,1/2)$ an asymptotic profile of the solution to \eqref{ike-100} is captured as 
\begin{equation}\label{ike-20}
\frac{e^{-t\vert\xi\vert^{2(1-\theta)}}}{\vert\xi\vert^{2\theta}}
\end{equation}
in \cite[Theorem 1.5]{IT}. In some sense, \eqref{ike-20} is similar to $\varphi_{1}(t,\xi)$ because of $\log(1+r^{2\theta}) \sim r^{2\theta}$ for small $r > 0$.}
\end{rem}
\begin{rem}{\rm A restriction $\theta \in (0,\displaystyle{\frac{1}{3}}]$ or $\theta \in (0,\displaystyle{\frac{5}{12}}]$ is just a technical condition, however, in the course of proof of Theorem \ref{maintheorem} one has frequently used the following fact
\begin{equation}\label{ike-11}
 \lim _{r \rightarrow +0} \frac{\log (1+r^{2\theta})}{r} = \infty.
\end{equation}
\eqref{ike-11} is also true in a more wider range $\theta \in (0,\frac{1}{2})$. So, reconsidering, the case of $\theta \in (\frac{1}{3}, \frac{1}{2})$ for $n = 1$ or $\theta \in (\frac{5}{12}, \frac{1}{2})$ for $n \geq 2$ is still open. }
\end{rem}
\begin{rem}{\rm The condition $u_{1} \in L^{2}({\bf R}^{n})$ in Theorem  \ref{maintheorem} is used to make sure the unique existence of the mild solution $u(t,x)$. However, it does not affect directly on the $L^{2}$-estimate of the solution, even in the high-frequency estimates although the estimate \eqref{eq3.50} in the high frequency zone can be easily estimated in terms of $ || u_1 || $ instead of $|| u_1 || _1$.}
\end{rem}
As an application of Theorem \ref{maintheorem} one can derive the following sharp decay estimates, which imply the optimal decay rates of the $L^{2}$-norm of the solution to problem \eqref{eqn}-\eqref{initial}.
\begin{theo}\label{main-teo} 
Let $n = 1$ with $0< \theta < \frac{1}{4}$ and $n\geq 2$ with $0< \theta \leq \frac{5}{12}$.  For $ u_1 \in L^{1,2\theta}({\bf R}^n)\cap L^{2}({\bf R}^n)$, it holds that
$$  K_1 \vert P_1\vert t^{-\frac{n-4\theta}{4(1-\theta)}} \leq \|u(t, \cdot)\| \leq K_2 (\vert P_1\vert+\|u_1\|_{L^{1,2\theta}}) \left ( t^{-\frac{n-4\theta}{4(1-\theta)}} + \frac{1}{\sqrt{\theta}} t^{-\frac{n-4\theta}{4\theta}}\right ) , \quad t\gg 1$$
with some constant $K_{1}, K_{2} > 0$ depending only on $n$ and $\theta$, where $u(t,x)$ is a unique solution to problem {\rm (1.1)-(1.2)} with $u_{0} = 0$. \\
\end{theo}

\begin{rem}{\rm A similar $L^{p}$-$L^{q}$ type "decay" estimates only from above has been already studied precisely in \cite{DE-2} and \cite[Corollary 2.2]{DE-3} to the solution of the equation \eqref{ike-100} for $n = 1$ and $0 < \theta < 1/4$, or $n \geq 2$ and $0 < \theta < 1/2$. The lower bound itself in Theorem \ref{main-teo} seems new.}
\end{rem}
\begin{rem}{\rm  As a result of Theorem \ref{main-teo}, one can observe that $\Vert u(t,\cdot)\Vert \sim t^{-\frac{n-4\theta}{4(1-\theta)}}$ ($t \to \infty$). Thus, as for an ultimate situation when $\theta \to 0^+$ formally, the optimal decay order will approach $t^{-\frac{n}{4}}$, which is the Gauss kernel. This is quite natural because in the case when $\theta = 0$, the equation corresponds to the frequently studied damped wave equation. In this sense, all results in this paper reflect a diffusive aspect of the equation (1.1) with small $\theta$. This property is quite different from those studied in \cite{CDI} for large $\theta \geq \frac{1}{2}$. In \cite{CDI}, a wave like property is captured.}
\end{rem}
\begin{rem}{\rm  The optimal decay order $t^{-\frac{n-4\theta}{4(1-\theta)}}$ obtained in Theorem \ref{main-teo} has a close relation to that studied in \cite{DE-3}, \cite{CLI} and \cite[(1.13) with $l=k= 0$]{IT} for the equation \eqref{ike-100}. In particular, the (almost) optimal decay rate of the "energy" and $L^{2}$-norm of the solutions are studied by developing a new energy method in the Fourier space in \cite{CLI}. So, the structure of the equation (1.1) is quite similar to \eqref{ike-100} with $\theta \in (0,1/2)$.}
\end{rem}

Contrary to the decay results as in Theorem \ref{main-teo}, one can observe the following surprising property, which shows infinite time blowup results of the solution to problem (1.1)-(1.2) in the one dimensional case. We believe this is the first discovery in the damped wave equation community. In \cite{DR} and \cite{DE-2}, when they apply the decay estimates of the solution for the equation \eqref{ike-100} to the nonlinear problems, they necessarily avoid to treat the case of $n = 1$, and $1/4 \leq \theta < 1/2$. The following crucial result makes their mechanism clear because of $\log(1+r^{2\theta}) \sim r^{2\theta}$ for small $r > 0$. 
\begin{theo}\label{main-teo-2} 
Let $n = 1$ with $\frac{1}{4} \leq \theta \leq \frac{1}{3}$.  For $ u_1 \in L^{1,2\theta}({\bf R})\cap L^{2}({\bf R})$, there exists positive constants $ K_1, K_2 $, which depend only on $\theta$, such that
\begin{align}\label{1.16}
K_1 |P_1| t^{\frac{4\theta -1 }{4\theta }} \leq \| u(t, \cdot ) \| \leq K_2 \big(  \frac{1}{\sqrt{4\theta -1}} |P_1| + \| u_1\|_{L^{1,2\theta}}  \big) t^{\frac{4\theta -1 }{4\theta}}, \quad t\gg1
\end{align}
for $ \frac{1}{4}< \theta \leq \frac{1}{3}  $ and 
\begin{align}\label{1.17}
K_1 |P_1| \sqrt{\log t}  \leq \| u(t, \cdot ) \| \leq K_2 \big( |P_1| + \| u_1\|_{L^{1,2\theta}}  \big)\sqrt{ \log t}, \quad t\gg1
\end{align}
for the case $ \theta=\frac{1}{4} $. 
\end{theo}
\begin{rem}{\rm We find that the number $\theta^{*} = \displaystyle{\frac{1}{4}}$ is critical in the one dimensional case because $\theta^{*}$ divides the structure of the corresponding solution $u(t,x)$ into two parts: one is decay property for $0 < \theta < \theta^{*}$, while the other is the infinite time blowup results in the case of $\theta^{*} \leq \theta < \frac{1}{2}$.
Moreover, we note that  in Theorem \eqref{main-teo-2} there is not a contradiction between the estimate  \eqref{1.16}  when $\theta \rightarrow (1/4)^{+}$   and the estimate \eqref{1.17} for $\theta=1/4$ because of the singularity $\sqrt{\frac{1}{4\theta-1}}$ at $\theta = 1/4$.
}
\end{rem}



{\bf Notation.} {\small Throughout this paper, $\| \cdot\|_q$ stands for the usual $L^q({\bf R}^{n})$-norm. For simplicity of notation, in particular, we use $\| \cdot\|$ instead of $\| \cdot\|_2$. Furthermore, we denote $\Vert\cdot\Vert_{H^{l}}$ as the usual $H^{l}$-norm. We  also define a relation $f(t) \sim g(t)$ as $t \to \infty$ by: there exist constant $C_{j} > 0$ ($j = 1,2$) such that
\[C_{1}g(t) \leq f(t) \leq C_{2}g(t)\quad (t \gg 1).\] 

For $\Omega \subset {\bf R}^n$ we denote $f \approx g$  on $\Omega$, if and only if 
there are constants $K_1, K_2$ such that 
$$ K_1 f(y) \leq g(y) \leq K_2 f(y), \;\;\mbox{for all} \; y \in \Omega.$$

We also introduce the following weighted functional spaces for $\gamma >  0$:
\[L^{1,\gamma}({\bf R}^{n}) := \left\{f \in L^{1}({\bf R}^{n}) \; \bigm| \; \Vert f\Vert_{L^{1,\gamma}} := \int_{{\bf R}^{n}}(1+\vert x\vert^{\gamma})\vert f(x)\vert dx < +\infty\right\}.\]
Finally, we denote the surface area of the $n$-dimensional unit ball by $\omega_{n} := \displaystyle{\int_{\vert\omega\vert = 1}}d\omega$. 

}


\section{Basic preliminary results}
In this section we shall collect important lemmas to derive precise estimates of the several  quantities related to the solution to problem \eqref{eqn}-\eqref{initial}. These are already studied and developed in our previous works (see \cite{Log-damping, CDI}).

The following estimate for the function  
$$I_p(t)= \int_0^{1}(1+r^{2})^{-t}r^p dr$$ 
is a direct consequence of the cases $p \geq 0$ in Char\~ao-Ikehata \cite{Log-damping} and $-1<p<0$ in Char\~ao-D'Abbicco-Ikehata \cite{CDI}.

\begin{lem}\label{general-p}
 Let $p > -1$ be a real number. Then
$$I_p(t) \sim t^{-\frac{p+1}{2}}, \quad t \gg 1.$$
\end{lem}

In order to deal with the high frequency part of estimates, one defines a function again 
$$J_p(t)=\int_1^{\infty}(1+r^2)^{-t}r^p dr$$
for $p \in {\bf R}$.

Then the next lemma is important to get estimates on the zone of high frequency to the solutions of the problem \eqref{eqn}--\eqref{initial}. The proof  appears in Char\~ao-Ikehata \cite{Log-damping}.
\begin{lem}\label{infit}
\,Let $p \in {\bf R}$. Then it holds that 
$$J_p(t) \sim \dfrac{2^{-t}}{t-1}, \quad t \gg 1.$$
\end{lem}
\vspace{0.2cm}
For later use we prepare the following simple lemma, which implies the exponential decay estimates of the middle frequency part.
\begin{lem}\label{intermid}\,Let $p \in {\bf R}$, and $\eta \in (0,1]$. Then there is a constant $C > 0$ such that 
$$\int_{\eta}^{1}(1+r^{2})^{-t}r^{p}dr \leq C(1+\eta^{2})^{-t}, \quad t \geq 0.$$
\end{lem}

\begin{rem}
{\rm We note that the proof of Lemma \ref{general-p} is done by using simple differential calculus and the theory from hypergeometric functions (see Watson \cite{W}). These are already developed in \cite{Log-damping} and \cite{CDI}.}

\end{rem}



\begin{lem}\label{lemmahiperbolicsine}
There exists a constant $K>0 $ such that 
$$ \frac{\sinh x }{x} \leq K\cdot e^{x}  $$
for $ x >0 $. 
\end{lem}

We will need the  following decomposition for the Fourier transform of a function $f$   in  $ L^1({\bf R}^n)$  as follows   
\begin{equation} \label{decompo}
\hat{f}(\xi)=A_f(\xi)-iB_f(\xi)+P_f,
\end{equation}
for all $\xi \in {\bf R}^n$ , where 
\begin{itemize}
\item[$\bullet$] $A_f(\xi)=\dfrac{1}{(2\pi)^{n/2}}\displaystyle{\int_{{\bf R}^n}}{(\cos(x\cdot\xi)-1)f(x)}dx,$
\item[$\bullet$] $B_f(\xi)=\dfrac{1}{(2\pi)^{n/2}}\displaystyle{\int_{{\bf R}^n}}{\sin(x\cdot\xi)f(x)}dx,$
\item[$\bullet$] $P_f=\dfrac{1}{(2\pi)^{n/2}}\displaystyle{\int_{{\bf R}^n}}{f(x)}dx.$
\end{itemize}
\noindent
Then, the next lemma has been already prepared in \cite{I-04} (see Notation for the definition of $L^{1,\kappa}({\bf R}^{n})$).
\begin{lem}\label{lema2.6}
\begin{itemize}
\item[{\rm i)}] If\;  $f \in L^1({\bf R}^n)$  then for all $\xi \in {\bf R}^n$ it is true that  
$$|A_f(\xi)|\leq L\|f\|_{L^1} \quad \text{ and  } \quad |B_f(\xi)|\leq N\|f\|_{L^1}.$$
\item[\rm{ii)}] If \;$0<\kappa \leq 1$ and  $f \in  L^{1,\kappa}({\bf R}^n)$ then for all  $\xi \in {\bf R}^n$ it is true that 
$$|A_f(\xi)|\leq K|\xi|^\kappa\|f\|_{L^{1,\kappa}} \quad \text{ and  } \quad |B_f(\xi)|\leq M|\xi|^\kappa\|f\|_{L^{1,\kappa}}$$
\end{itemize}
\noindent with  $L$, $N$, $K$ and  $M$ positive  constants  depending only on  the dimension  $n$  or $n$ and $\kappa$.\\
\end{lem}

\begin{lem}\label{lema2-2theta}
Let $0\leq \theta < 1$ and $q>-1$. Then 
$$  \int_0^1 (1+r^{2-2\theta})^{-t}r^q dr \sim \dfrac{1}{1-\theta}  t^{-\frac{q+1}{2(1-\theta)}}, \quad t\gg 1 .$$
In particular, for $0\leq \theta \leq 1/2$ and $q>-1$ it holds that $$  \int_0^1 (1+r^{2-2\theta})^{-t}r^q dr \sim t^{-\frac{q+1}{2(1-\theta)}}, \quad t\gg 1 .$$
\end{lem}
{\it Proof.}\ Let $ s = r^{1-\theta} $. Then
\begin{equation*}
\int_0^1 (1+r^{2-2\theta})^{-t}r^q dr= \frac{1}{1-\theta}  \int _0^1 (1+s^{2})^{-t} s^{\frac{q+\theta}{1-\theta}} ds
\end{equation*}
Since $0\leq \theta < 1$ and $q>-1 $, we have $\frac{q+\theta}{1-\theta}>-1$. Thus, we can apply the Lemma \ref{general-p} to obtain the result. 
\hfill
$\Box$
\\
\begin{rem}\label{estimativaporbaixoIp} {\rm Actually, for $ \eta >0  $, $0\leq \theta  \leq 1/2 $ and $q>-1$, it holds that 
$$\int_0^{\eta} (1+r^{2-2\theta})^{-t}r^q dr \geq C t^{-\frac{q+1}{2(1-\theta)}}, \quad t\gg 1 $$
for some constant $C>0$ depending on each $\eta > 0$.}
\end{rem} 
Indeed, it suffices to check the case of $0 < \eta < 1$.  In this case, one notices
\[\int_0^\eta (1+r^{2-2\theta})^{-t}r^q dr = \int_0^1 (1+r^{2-2\theta})^{-t}r^q dr - \int_\eta^1 (1+r^{2-2\theta})^{-t}r^q dr,\]
and one has
\[\int_\eta^1 (1+r^{2-2\theta})^{-t}r^q dr \leq \frac{1}{1+q}(1-\eta^{q+1})(1+\eta^{2-2\theta})^{-t}.\]
Since the last term implies the exponential decay, the desired estimate can be derived soon via Lemma \ref{lema2-2theta}.
\hfill
$\Box$


\begin{lem}\label{lema2theta}
Let $ \theta >0 $ and $q>-1$. Then
$$  \int_0^1 (1+r^{2\theta})^{-t}r^q dr \sim \frac{1}{\theta}t^{-\frac{q+1}{2\theta}}, \quad t\gg1 .$$
\end{lem}
{\it Proof.}\
We consider the change of variable $ s=r^{\theta} $. Then 
\begin{equation*}
 \int_0^1 (1+r^{2\theta})^{-t}r^q dr = \frac{1}{\theta}  \int_0^1 (1+s^{2})^{-t}s^{\frac{q+1-\theta}{\theta}} ds
\end{equation*}
for $t \geq 0.$ Finally, $ \frac{q+1-\theta}{\theta} >-1$ because of $ q>-1 $. From Lemma \ref{general-p} the desired result follows. 
\hfill
$\Box$
\\

\begin{lem}\label{lema2thetainfit}
Let $ \theta >0 $ and $q \in {\bf R}$. Then
$$  \int_{1}^{\infty} (1+r^{2\theta})^{-t}r^q dr \sim \frac{1}{\theta}  \frac{2^{-t}}{t-1}, \quad t \gg 1 .$$
\end{lem}

\begin{lem}\label{lema2-2thetainfit}
Let $0\leq \theta < 1$ and $q\in {\bf R}$. Then 
$$  \int_1^{\infty} (1+r^{2-2\theta})^{-t}r^q dr \sim  \frac{2^{-t}}{t-1}, \quad t \gg 1.$$
\end{lem}

{\it Proof of Lemmas \ref{lema2thetainfit} and \ref{lema2-2thetainfit}.}\, From lemma \ref{infit} and the change of variables as in Lemmas \ref{lema2theta} and \ref{lema2-2theta}, the result now follows. 
\hfill
$\Box$

\section{Asymptotic profile}

As mentioned in the introduction, our main interest in this work is to investigate the problem \eqref{eqn}-\eqref{initial} for the case of $0<\theta<1/2$ in order to compensate the research in \cite{CDI} studying the case of $\theta>1/2$.

The associated problem to \eqref{eqn}-\eqref{initial}  in Fourier space is the following
\begin{align}
& \hat{u}_{tt} + \log (1+|\xi|^{2\theta})\hat{u}_t +|\xi|^2 \hat{u} = 0,\ \ \ t>0,\quad \xi \in {\bf R}^{n},\label{eqnfourier}\\
& \hat{u}(0,\xi)= 0, \quad  \hat{u}_{t}(0,\xi)= \hat{u}_1(\xi),\ \ \ \xi \in{\bf R}^{n} ,\label{initialfourier}
\end{align}
where the associated  characteristic polynomial is 
$$ \lambda^2  + \log (1+|\xi|^{2\theta}) \lambda + |\xi|^2 = 0. $$
The characteristics roots are expressed as 
\begin{equation}\label{defilambda}
\lambda _{\pm} = \frac{- \log (1+|\xi|^{2\theta}) \pm \sqrt{\log ^2(1+|\xi|^{2\theta}) - 4 |\xi|^2 }}{2}, \quad \xi \in {\bf R}^n . 
\end{equation}

\begin{lem}\label{lemadelta}
There exists $ \delta =\delta(\theta) $, $0<\delta<1$ such that 
\begin{align}
&\log ^2(1+|\xi|^{2\theta}) - 4 |\xi|^2 \geq 0  \text{ for } |\xi| \leq \delta, \label{definicao1delta}\\ 
& \log ^2(1+|\xi|^{2\theta}) - 4 |\xi|^2 < 0  \text{ for }  |\xi| > \delta. \label{definicao2delta}
\end{align}
\end{lem}
{\it Proof.}\
Working with $r=|\xi| $,  we first  observe that $ \log(1+r^2) < 2r  $ for all $r> 0 $. Also, $r^{2\theta} \leq r^2$  for $ r \geq 1 $, since $\theta \leq 1$. Therefore, in the case $\theta < \displaystyle{\frac{1}{2}}$, one has
\begin{equation}\label{deltaleq1}
\log (1+r^{2\theta}) \leq \log (1+r^2) < 2r
\end{equation}
for all $ r \geq 1 $.  Thus we may conclude that the function $ f(r):= \log (1+r^{2\theta}) -  2r $ is negative for all $ r \geq 1 $. 
However, the similar phenomena does not happen near the origin. In fact, we first notice that 
$$ \lim _{r \rightarrow +0} \frac{\log (1+r^{2\theta})}{r} = \infty, $$
for  $ \theta \in (0,\frac{1}{2}) $. Therefore, there exists $r_{0} = r_{0}(\theta) <1$ such that 
\begin{equation*}
\frac{\log (1+r^{2\theta})}{r} >2
\end{equation*}
for all $ r\in {\bf R}^n$ satisfying $ 0<r< r_{0}$.  Then,  $f(r) = \log (1+r^{2\theta}) - 2r \geq 0 $ for $0\leq r<r_{0}$. Furthermore, for $ r \geq 0 $ one can get  
\begin{align*}
f''(r) &= \frac{2\theta r^{2\theta-2} \left [ (2\theta -1)(1+r^{2\theta})-2\theta r^{2\theta} \right ]  }{(1+r^{2\theta})^2} = \frac{2\theta r^{2\theta-2} \left [ 2\theta -1 - r^{2\theta} \right ]  }{(1+r^{2\theta})^2}. 
\end{align*}
Since  $ 0 <\theta<\displaystyle{\frac{1}{2}}$, the function $ f: [0, \infty) \rightarrow {\bf R}$ satisfies $f''(r) < 0$. Due to $f(0)=0$ and \eqref{deltaleq1} one can conclude that there exists a unique number $\delta=\delta (\theta) $, $ 0 < \delta<1 $ such that $ f(r) \geq 0 $ for all $ 0\leq r \leq \delta $ and $ f(r) \leq 0 $ for all $ r\geq \delta $. Finally, one can write 
$$ \log ^2(1+|\xi|^{2\theta}) - 4 |\xi|^2 = f(|\xi|) \left ( \log (1+|\xi|^{2\theta}) + 2 |\xi| \right ).$$ 
Therefore, using the properties of the function $f(r)=f(|\xi|)$ one can obtain the desired statement. 
\hfill
$\Box$
\\
By Lemma \ref{lemadelta}, we see that that the characteristics roots \eqref{defilambda} are real-valued for $ |\xi|\leq \delta $, and complex-valued for $ |\xi| > \delta $. This is a crucial different point from that observed in the case of $1/2 \leq \theta$. 

%

\subsection{Estimates on the region $ |\xi| \leq \delta $}

First part of this section we analyze the behavior of the characteristics roots near the origin $\xi = 0$. To do that we  need some remarks  and lemmas. 

\begin{rem}\label{rem3.1} {\rm For $q\geq 0 $ it is easy to check the inequality $\frac{1}{2} r^q \leq \log (1+r^q) \leq  r^q $ for $r \in [0,1]$.\\
In particular, for  $ 0<  \theta < \frac{1}{2}$ we have 
\begin{align}
& \frac{1}{2} |\xi|^{2\theta} \leq \log(1+|\xi|^{2\theta}) \leq \frac{3}{2}|\xi|^{2\theta},\label{defieta1} \\
& \frac{1}{2} |\xi|^2 \leq \log(1+|\xi|^2) \leq \frac{3}{2}|\xi|^2,\\
& \frac{1}{2} |\xi|^{2-2\theta} \leq \log(1+|\xi|^{2-2\theta}) \leq \frac{3}{2}|\xi|^{2-2\theta} \label{defieta3}
\end{align}
for $ |\xi| \leq 1 $.}
\end{rem}


We note that for $0\leq \theta<1/2$ it holds that
$$ \lim _{r \rightarrow +0} \frac{r^{4-4\theta}}{r^2} = 0 .$$
Thus, there exists $ \delta _1=\delta_1(\theta)$, $ 0<\delta_1 <1$ that satisfies 
\begin{equation}\label{defidelta2}
\frac{|\xi|^{4-4\theta}}{|\xi|^2} \leq \frac{1}{25}
\end{equation}
whenever $ 0<|\xi|\leq \delta_1 $.  Moreover, one can choose $ \delta_1 \in (0,1)$ such as $\delta_{1} < \delta $. In fact, from \eqref{defidelta2} one has
$$ 25 |\xi|^2 \leq |\xi|^{4\theta}$$
for $ 0\leq |\xi| \leq  \delta_1 $. In this region it also holds  $ |\xi|^{4\theta} \leq 4\log^2(1+|\xi|^{2\theta}) $, due to \eqref{defieta1}. 
Thus 
\begin{equation} \label{25-4}
\log^2(1+|\xi|^{2\theta}) \geq \frac{25}{4}|\xi|^2 \geq \frac{16}{3}|\xi|^2  \geq 4 |\xi|^2 
\end{equation}
for $|\xi|\leq \delta_1. $  Comparing \eqref{25-4} with \eqref{definicao1delta}, we may conclude that $\delta_1< \delta$. From \eqref{25-4} we also obtain 
  $$ \log^2(1+|\xi|^{2\theta}) \geq \frac{16}{3}|\xi|^2  $$
 whenever $|\xi| \leq \delta_1$.\\
  
Now we define a new number: 
 \begin{equation}\label{etainfimo}
 \eta := \sup \{ \alpha >0; \frac{|\xi|^{4-4\theta}}{|\xi|^2} \leq \frac{1}{25} \text{ for } 0<| \xi|\leq \alpha  \} .
 \end{equation}
We note that $\eta$ is positive and is well defined, because the set $ \{ \alpha >0; \frac{|\xi|^{4-4\theta}}{|\xi|^2} \leq \frac{1}{25} \text{ for } | \xi|\leq \alpha  \}  $ is not empty  ($\delta_1$ is a member of this set) and is bounded from above. In fact, for example,  $1$ is an upper bound for this set, and $\eta<\delta<1$ with $\delta$ defined in Lemma \ref{lemadelta}. In particular, the following two properties are true for $|\xi|\leq \eta$:
\begin{align}
& \frac{3}{4} \log^2(1+|\xi|^{2\theta}) \geq 4|\xi|^2, \label{defieta4}\\
& 25 |\xi|^{4 -4\theta} \leq |\xi|^2 . \label{defieta5}
\end{align}

\begin{lem}\label{lemmaequivalencias}  Let $ \eta$  be  the number defined by \eqref{etainfimo}. Then, for  $|\xi| \leq \eta $ it holds that
\begin{itemize}
\item[{\rm (i).}]\,\,$ \lambda_+ - \lambda_- \approx \log (1+|\xi|^{2\theta});$
\item[{\rm (ii).}]\,$ \lambda _+ \approx -\log (1+|\xi|^{2-2\theta}) \approx -|\xi|^{2-2\theta} ; $
\item[{\rm (iii).}]$\lambda_- \approx -\log (1+|\xi|^{2\theta}) . $
\end{itemize}
\end{lem}
{\it Proof.}\\
{\rm (i)}\,The upper estimate is simple  because for $| \xi| \leq \eta < \delta $ it holds that 
$$ \lambda_+ - \lambda _-= \sqrt{\log ^2(1+|\xi|^{2\theta}) - 4 |\xi|^2} \leq \sqrt{\log ^2(1+|\xi|^{2\theta})} =\log (1+|\xi|^{2\theta}) . $$
On the other hand, by \eqref{defieta4} we have
 $$ \frac{1}{4} \log^2 (1+|\xi|^{2\theta}) \leq \log^2 (1+|\xi|^{2\theta}) - 4 |\xi|^2, \quad |\xi|\leq \eta  .$$
 For this reason, in the zone $|\xi|\leq \eta $ it holds that 
 \begin{equation*}
  \frac{1}{2} \log (1+|\xi|^{2\theta}) \leq \sqrt{\log^2 (1+|\xi|^{2\theta}) - 4 |\xi|^2} .
 \end{equation*}

{\rm (ii).}\,The inequality \eqref{defieta5} provides us to get
 \begin{equation}\label{eq3.13}
 0\geq 25 |\xi|^{4-4\theta} - 5 |\xi|^2 +4|\xi|^2 = 25 |\xi|^{4-4\theta} - 5 |\xi|^{2\theta}|\xi|^{2-2\theta} +4|\xi|^2 .
 \end{equation}
 The lower inequality in  \eqref{defieta1} implies  that $ -10 \log(1+|\xi|^{2\theta}) \leq -5 |\xi|^{2\theta} $ for $ |\xi|\leq 1$ and in particular for  $ |\xi|\leq \eta $. By combining this fact with \eqref{eq3.13}, we obtain 
 \begin{equation*}
 25 |\xi|^{4-4\theta} - 10 \log (1+|\xi|^{2\theta}) |\xi|^{2-2\theta} +4|\xi|^2  \leq 0, \quad |\xi|\leq \eta .
 \end{equation*}
Adding $\log^2 (1+|\xi|^{2\theta})$ on both sides we may obtain 
$$ \left ( \log (1+|\xi|^{2\theta}) - 5 |\xi|^{2-2\theta} \right )^2  = \log^2 (1+|\xi|^{2\theta}) - 10 \log (1+|\xi|^{2\theta}) |\xi|^{2-2\theta} + 25 |\xi|^{4-4\theta} \leq  \log^2 (1+|\xi|^{2\theta}) - 4|\xi|^2. $$
Hence, for $|\xi|\leq \eta $,  $  \log (1+|\xi|^{2\theta}) - 5 |\xi|^{2-2\theta} \leq \sqrt{\log^2 (1+|\xi|^{2\theta}) - 4|\xi|^2} $ and
$$ -\frac{5}{2} |\xi|^{2-2\theta} \leq \frac{-\log (1+|\xi|^{2\theta}) + \sqrt{\log^2 (1+|\xi|^{2\theta}) - 4|\xi|^2} }{2} = \lambda _+ .$$
Furthermore, we also concludes that
\begin{equation}\label{upper-ii} -5 \log(1+|\xi|^{2-2\theta}) \leq  -\frac{5}{2} |\xi|^{2-2\theta} \leq  \lambda _+
\end{equation}
on the zone $|\xi| \leq \eta $, due to  \eqref{defieta3}.

In order to prove the upper estimate part of {\rm (ii)} we first observe that 
\begin{equation*}
0 \leq |\xi|^2 +|\xi|^{4-4\theta} = 4 |\xi|^2 +|\xi|^{4-4\theta} - 3 |\xi|^{2\theta} |\xi|^{2-2\theta}.
\end{equation*}
In the zone  $|\xi|\leq \eta $ it holds that $ -3 |\xi|^{2\theta} \leq -2 \log (1+|\xi|^{2\theta}) $ by \eqref{defieta1}, which implies that
$$ -3 |\xi|^{2\theta} |\xi|^{2-2\theta} \leq -2 \log (1+|\xi|^{2\theta}) |\xi|^{2-2\theta}.$$ 

By using the inequality just above we may obtain that 
\begin{equation}\label{316}
0 \leq 4|\xi|^2 +|\xi|^{4-4\theta} -2 \log (1+|\xi|^{2\theta}) |\xi|^{2-2\theta}. 
\end{equation}
We add $ \log ^2(1+|\xi|^{2\theta})  $ in both side of \eqref{316} in order to get the following estimate:
\begin{equation*}
\log^2(1+|\xi|^{2\theta}) -  4|\xi|^2 \leq \log^2(1+|\xi|^{2\theta})   -2 \log (1+|\xi|^{2\theta}) |\xi|^{2-2\theta} +|\xi|^{4-4\theta} = \left ( \log (1+|\xi|^{2\theta}) - |\xi|^{2-2\theta}  \right )^2.
\end{equation*}
This implies
\begin{equation}\label{eq317}
\lambda _+=\frac{- \log (1+|\xi|^{2\theta})+ \sqrt{\log^2(1+|\xi|^{2\theta}) -  4|\xi|^2 }}{2} \leq    - \frac{1}{2}|\xi|^{2-2\theta}.
\end{equation}
When one derives \eqref{eq317}, one must check the fact that $\log (1+|\xi|^{2\theta}) - |\xi|^{2-2\theta} \geq 0$ on $\vert\xi\vert \leq \eta$. Indeed, this can be easily observed by a combination of \eqref{defieta4} and \eqref{defieta5}.

Now, by combining inequalities \eqref{eq317} and \eqref{defieta3} one obtain 
$$ \lambda _+ \leq    - \frac{1}{2}|\xi|^{2-2\theta} \leq -\frac{1}{3} \log(1+ |\xi|^{2-2\theta} ) $$
because of $|\xi|\leq \eta$. The inequalities just above and \eqref{upper-ii} imply the desired statement of item {\rm (ii)}. \\

{\rm (iii).}\,In the course of the proof of item {\rm (i)} in the region $|\xi|\leq \eta$, we also have
$$ - \log (1+|\xi|^{2\theta}) \leq - \sqrt{\log^2 (1+|\xi|^{2\theta}) - 4 |\xi|^2} \leq - \frac{1}{2} \log (1+|\xi|^{2\theta}).$$
Therefore, one can easily conclude that 
$$ - \log (1+|\xi|^{2\theta}) \leq \frac{ - \log (1+|\xi|^{2\theta}) - \sqrt{\log^2 (1+|\xi|^{2\theta}) - 4 |\xi|^2}}{2} \leq - \frac{3}{4} \log (1+|\xi|^{2\theta}), \quad   |\xi| \leq \eta. $$
This implies the desired statement of item {\rm (iii)}.
\hfill
$\Box$

\subsubsection{Estimates on the low-frequency zone	$|\xi| \leq \eta^3$}
Throughout this paper we assume the initial amplitude $u_{0}$ satisfies $u_0=0$,  without loss of generality in order to investigate the asymptotic profiles of solutions.
  
We first remember the number $\eta \in (0,\delta)$ defined in \eqref{defieta1}-\eqref{defieta5}. Also, since $ 0<\eta < 1 $, we have $ \eta ^3 < \eta $. In the zone of low frequency $\vert\xi\vert \leq \eta^3<\eta$, the characteristics roots $\lambda_{\pm} $ are real, and the solution of \eqref{eqnfourier}-\eqref{initialfourier} is explicitly given by 
\begin{equation}\label{solutioneta}
\hat{u}(t,\xi)= \frac{e^{t\lambda_+}-e^{t\lambda_-}}{\lambda_+-\lambda_-}\hat{u}_1(\xi).
\end{equation}

The purpose in this section is to get an asymptotic profile to the solution $\hat{u}(t,\xi)$, and in order to do that we need to obtain useful estimates. For this reason, we defined a function $ g:[0,\delta] \rightarrow {\bf R}$ inspired by an idea from \cite{FIM-1}, as follows. A discovery of this function $g(s)$ is one of decisive points in our proof.
\begin{align}\label{defig}
g(s):= \left\{\begin{matrix}
1+ \sqrt{1- \frac{4 s^6}{\log^2(1+s^{6\theta})}} & \text{ if } 0<s\leq \delta \\ 
2 & \text{ if } s=0.
\end{matrix}\right.
\end{align}

Note that for $0 < \theta < 1/2$,
\[\lim_{s \to 0^+}\frac{s^6}{\log^2(1+s^{6\theta})} = 0.\]

\begin{rem}\label{RemarkA1} {\rm Let $ t>0$ and $ \xi \in {\bf R}^n  $, $ 0< |\xi|\leq \eta$,  be fixed.  We recall that $\eta < \delta<1$. Let us consider the function $h(s)$ defined on $[0,\eta]$ as follows:
\begin{align*}
h(s) := e^{-\frac{t\log(1+\vert\xi\vert^{2\theta})}{2}g(s)}.
\end{align*}
We see that $h(s)$ is differentiable on $(0,\eta)$. Then, it should be noted that one can apply the mean value theorem in the interval $[0, s]$ for each $s \in (0,\eta]$ to get
\begin{align}\label{A3.25}
\frac{h(s)-h(0)}{s} = \frac{e^{ -\frac{t\log(1+|\xi|^{2\theta})}{2} g(s)} - e^{-\frac{t\log(1+|\xi|^{2\theta})}{2} g(0)}}{s} = 
 -\frac{t\log(1+|\xi|^{2\theta})}{2}  e^{ -\frac{t\log(1+|\xi|^{2\theta})}{2} g(\alpha s  )} g'(\alpha s)
\end{align}
with some  $\alpha=\alpha(s, t, |\xi|) \in (0,1 )$.}
\end{rem} 

We observe that on the low frequency zone $ 0\leq |\xi| \leq \eta ^3 $ it holds that
$$ \lambda _- = -\frac{\log (1+|\xi|^{2\theta})}{2} g(\sqrt[3]{|\xi|}).$$
By applying \eqref{A3.25} for $t>0$ and  $ s=\sqrt[3]{|\xi|}  $, $0< | \xi |  \leq \eta ^3$,  we have 
\begin{equation}\label{decomp1}
e^{t\lambda_{-}}=e^{-t \log(1+|\xi|^{2\theta})} -\frac{t}{2} \log(1+|\xi|^{2\theta})\sqrt[3]{|\xi|}  e^{ -\frac{t\log(1+|\xi|^{2\theta})}{2} g(\alpha \sqrt[3]{|\xi|}  )} g'(\alpha \sqrt[3]{|\xi|})
\end{equation}
with $\alpha := \alpha(s, t, \vert\xi\vert) = \alpha (t,|\xi|) \in (0,1) $.
\vspace{0.2cm}

From the Chill-Haraux \cite{C-H-01} idea, we also observe that 
$$ \lambda _+=-\frac{\lambda_+^2+|\xi|^2}{\log(1+|\xi|^{2\theta})},$$
so that one has 
\begin{align}\label{decomp2}
e^{t\lambda _+} =  e^{-\frac{|\xi|^2}{\log(1+|\xi|^{2\theta})}t} e^{- \frac{\lambda_+^2}{\log(1+|\xi|^{2\theta})}t} .
\end{align}
On the other hand, because of \eqref{defilambda} we see that 
 $$ \frac{1}{\lambda_+-\lambda_-}=\frac{1}{\log(1+|\xi|^{2\theta})} +R(|\xi|) $$
where
\begin{equation}\label{defiR(xi)}
R(r)=  \frac{4 r^2}{ \log^3(1+r^{2\theta}) \sqrt{1 - \frac{4 r^2}{\log^2(1+r^{2\theta})}} \left ( 1+  \sqrt{1- \frac{4 r^2}{\log^2(1+r^{2\theta})}} \right ) } .
\end{equation}

By combining \eqref{decomp1}, \eqref{decomp2} and \eqref{defiR(xi)} with the decomposition of initial data
$$ \hat{u}_1 (\xi) = A_{u_{1}}(\xi)-iB_{u_{1}}(\xi) +P_{u_{1}} =:  A_1(\xi)-iB_1(\xi) + P_1$$
as in \eqref{decompo}, we can write the solution of \eqref{eqnfourier}-\eqref{initialfourier} given by \eqref{solutioneta}, for $|\xi| \leq \eta^3$, as follows
\begin{align}\label{solutiondecomplow}
\hat{u}(t,\xi) &= \frac{e^{-\frac{|\xi|^2}{\log(1+|\xi|^{2\theta})}t}}{\log(1+|\xi|^{2\theta})} P_1  - \frac{e^{-t \log(1+|\xi|^{2\theta})}}{\log(1+|\xi|^{2\theta})} P_1 + R(|\xi|) e^{-\frac{|\xi|^2}{\log(1+|\xi|^{2\theta})}t} \hat{u}_1(\xi) - R(|\xi|) e^{-t \log(1+|\xi|^{2\theta})} \hat{u}_1(\xi) \nonumber \\
& + \frac{e^{-\frac{|\xi|^2}{\log(1+|\xi|^{2\theta})}t}}{\log(1+|\xi|^{2\theta})} (A_1(\xi) -i B_1(\xi)) + e^{-\frac{|\xi|^2}{\log(1+|\xi|^{2\theta})}t} \frac{ e^{- \frac{\lambda_+^2}{\log(1+|\xi|^{2\theta})}t} -1 }{\lambda_+ -\lambda_-} \hat{u}_1(\xi) \nonumber \\
&- \frac{e^{-t \log(1+|\xi|^{2\theta})}}{\log(1+|\xi|^{2\theta})}(A_1(\xi)-iB_1(\xi)) + t \frac{\log(1+|\xi|^{2\theta}) \sqrt[3]{|\xi|}  }{2 (\lambda_+-\lambda _-)}  e^{ -\frac{t\log(1+|\xi|^{2\theta})}{2} g\left ( \alpha \sqrt[3]{|\xi|} \right ) } g'\left ( \alpha \sqrt[3]{|\xi|} \right ) \hat{u}_1(\xi).
\end{align}

Our main goal in this subsection is to introduce an asymptotic profile as $t \to +\infty$ of the solution $\hat{u}(t,\xi)$ in the low frequency region  $|\xi|  \leq \eta^3< \eta$  in the simple form
\begin{equation}\label{defivarphi1}
\varphi (t,\xi) := \frac{e^{-\frac{|\xi|^2}{\log(1+|\xi|^{2\theta})}t}}{\log (1+|\xi|^{2\theta})} P_1 - \frac{e^{-\log(1+|\xi|^{2\theta})t}}{\log(1+|\xi|^{2\theta})}P_1.
\end{equation}
Thus, we need to prove that 
$$ \| \hat{u}(t,\cdot ) - \varphi (t,\cdot) \| \rightarrow 0 , \quad t\rightarrow \infty$$
much faster than the   components  in the right hand side of \eqref{defivarphi1}. For this purpose, we consider the following six remainder functions
\begin{align*}
&F_1(t,\xi) = R(|\xi|) e^{-\frac{|\xi|^2}{\log(1+|\xi|^{2\theta})}t} \hat{u}_1(\xi) \\
&F_2(t,\xi) = - R(|\xi|) e^{-t \log(1+|\xi|^{2\theta})} \hat{u}_1(\xi) \\
&F_3(t,\xi) = \frac{e^{-\frac{|\xi|^2}{\log(1+|\xi|^{2\theta})}t}}{\log(1+|\xi|^{2\theta})} (A_1(\xi) -i B_1(\xi)) \\
&F_4(t,\xi) = e^{-\frac{|\xi|^2}{\log(1+|\xi|^{2\theta})}t} \frac{ e^{- \frac{\lambda_+^2}{\log(1+|\xi|^{2\theta})}t} -1 }{\lambda_+ -\lambda_-} \hat{u}_1(\xi)  \\
&F_5(t,\xi) = - \frac{e^{-t \log(1+|\xi|^{2\theta})}}{\log(1+|\xi|^{2\theta})}(A_1(\xi)-iB_1(\xi)) \\
&F_6(t,\xi) = t \frac{\log(1+|\xi|^{2\theta}) \sqrt[3]{|\xi|}  }{2 (\lambda_+-\lambda _-)}  e^{ -\frac{t\log(1+|\xi|^{2\theta})}{2} g\left ( \alpha \sqrt[3]{|\xi|} \right ) } g'\left ( \alpha \sqrt[3]{|\xi|} \right ) \hat{u}_1(\xi).
\end{align*}

From \eqref{solutiondecomplow} and \eqref{defivarphi1}, for $|\xi| \leq \eta$, we have  
\begin{equation*}
\hat{u}(t,\xi) - \varphi(t,\xi)   = \sum _{j=1}^{6}  F_j(t,\xi) .
\end{equation*}

In order to obtain decay rates in time to these functions we assume the additional condition on the initial data such that
$$ u_1 \in L^{1,2\theta}({\bf R}^{n}), \,\,\, 0<\theta<1/2.$$ 
To begin with, we estimate the function $F_{3}(t,\xi)$. Indeed, by using Lemma \ref{lemmaequivalencias}, and Lemma \ref{lema2.6} with $\kappa := \theta \in (0,1/2)$ (this is our crucial idea), and the inequality \eqref{defieta1} and \eqref{defieta3} one can estimate 
\begin{align}\label{estimativaF_3}
\int _{|\xi|\leq \eta^3<1} |F_3(t,\xi)|^2d\xi &= \int _{|\xi|\leq \eta^3}\frac{e^{-\frac{2|\xi|^2}{\log(1+|\xi|^{2\theta})}t}}{\log^2 (1+|\xi|^{2\theta})} \vert A_1(\xi)-iB_1(\xi)\vert^2d\xi \nonumber \\
&\leq \int _{|\xi|\leq \eta^3}\frac{e^{-\frac{2t}{3} \log(1+|\xi|^{2-2\theta}) }}{\log^2 (1+|\xi|^{2\theta})} \vert A_1(\xi)-iB_1(\xi)\vert^2d\xi \nonumber \\
&\leq (M+K)^2  \|u_1\|_{L^{1,2\theta }}^2 \int _{|\xi|\leq \eta^3}\frac{e^{-\frac{2t}{3}\log(1+|\xi|^{2-2\theta}) }}{\log^2 (1+|\xi|^{2\theta})} |\xi|^{4\theta} d\xi \nonumber \\
&\leq C (M+K)^2  \|u_1\|_{L^{1,2\theta }}^2 \int _{|\xi|\leq \eta^3} (1+|\xi|^{2-2\theta})^{-\frac{2t}{3}}  d\xi \nonumber \\
&= C (M+K)^2 \omega_n \|u_1\|_{L^{1,2\theta }}^2 \int _0^{\eta^3} (1+r^{2-2\theta})^{-\frac{2t}{3}} r^{n-1} dr \nonumber \\
&\leq  C \|u_1\|_{L^{1,2\theta }}^2 t^{-\frac{n}{2(1-\theta)}}, \quad t \gg 1.
\end{align}
with a generous constant $C > 0$ depending only on $n$, and the last inequality is due to Lemma \ref{lema2-2theta}, where
one has just used the fact that
\[\lim_{\sigma \to +0}\frac{\sigma}{\log(1+\sigma)} = 1.\]

Similarly, we can also estimate 
\begin{align}\label{estimativaF_5}
\int_{|\xi|\leq \eta^3} |F_5(t,\xi)|^2d\xi &=\int_{|\xi|\leq\eta^3} \frac{e^{-2\log(1+|\xi|^{2\theta})t}}{\log^2(1+|\xi|^{2\theta})} \vert A_1(\xi)-iB_1(\xi)\vert^2 d\xi \nonumber \\
&= \int_{|\xi|\leq\eta^3} \frac{(1+|\xi|^{2\theta})^{-2t}}{\log^2(1+|\xi|^{2\theta})} \vert A_1(\xi)-iB_1(\xi)\vert^2 d\xi \nonumber \\
&\leq C (K+M)^2 \|u_1\|_{L^{1,2\theta}}^2 \int_{|\xi|\leq\eta^3} (1+|\xi|^{2\theta})^{-2t} d\xi \nonumber \\
&= C \omega _n (K+M)^2 \|u_1\|_{L^{1,2\theta}}^2 \int_0^{\eta^3} (1+r^{2\theta})^{-2t} r^{n-1} \nonumber dr \\
&\leq \frac{C}{\theta} \|u_1\|_{L^{1,2\theta}}^2 t^{-\frac{n}{2\theta}},
 \quad t \gg 1,
\end{align}
where the last inequality is due to Lemma \ref{lema2theta}.

On the next  estimates to the functions $F_j(t,\xi)$ we also rely on Lemma \ref{lema2-2theta} or Lemma \ref{lema2theta}.\\
In order to estimate $F_4(t,\xi)$ we use the fact that $ |e^{-a} -1 | \leq a $ for all $a \geq 0$. Then, Lemma \ref{lemmaequivalencias} and inequality 
\eqref{defieta3}
imply the existence of a constant $C>0$ such  that 
\begin{align}\label{F4rate}
\int_{|\xi| \leq \eta^3} |F_4(t,\xi)|^2 d\xi &= \int_{|\xi| \leq \eta^3} \left ( \frac{   e^{\frac{-\lambda_+^2}{\log(1+|\xi|^{2\theta})}t}-1  }{\lambda_+-\lambda _-}  \right )^2 e^{-\frac{2|\xi|^2}{\log(1+|\xi|^{2\theta})}t} |\hat{u}_1(\xi)|^2 d\xi \nonumber \\
&\leq t^2 \| u_1\|_{1}^2 \int_{|\xi| \leq \eta^3}  \frac{\lambda_+^4}{\log^2(1+|\xi|^{2\theta})} \frac{e^{-\frac{2|\xi|^2}{\log(1+|\xi|^{2\theta})}t}}{(\lambda_+-\lambda _-)^2} d\xi \nonumber \\
&\leq Ct^2 \| u_1\|_{1}^2 \int_{|\xi| \leq \eta^3}  \frac{|\xi|^{8-8\theta} }{\log^4(1+|\xi|^{2\theta})} e^{-\frac{2t}{3}\log(1+|\xi|^{2-2\theta})} d\xi \nonumber \\
&\leq Ct^2 \| u_1\|_{1}^2 \int_{|\xi| \leq \eta^3}  (1+|\xi|^{2-2\theta})^{-\frac{2t}{3}}  |\xi|^{8-16\theta} d\xi \nonumber \\
&= C\omega_n t^2 \| u_1\|_{1}^2 \int_{0}^{\eta^3}  (1+r^{2-2\theta})^{-\frac{2t}{3}} r^{7-16\theta +n} dr \nonumber \\
&\leq C  \| u_1\|_{1}^2 t^2 t^{-\frac{8-16\theta+n}{2(1-\theta)}} \nonumber \\
&= C  \| u_1\|_{1}^2  t^{-\frac{4-12\theta+n}{2(1-\theta)}}, \quad t \gg 1.
\end{align}

\begin{rem} {\rm Note that in the above estimate \eqref{F4rate} to apply Lemma \ref{lema2-2theta} it is necessary to check $7-16\theta+n >-1$, but this holds for $0\leq \theta <1/2$. Moreover, according to our computations above, we have to prove that all $L^2$-norm of the six functions $F_1(t,\xi), \cdots, F_6(t,\xi)$ decay to zero in time.  However, to get such decay estimates in \eqref{F4rate}, we need additional restriction such that $0\leq \theta <\frac{5}{12} < \frac{1}{2}$ in the case $n=1$. For $n \geq 2$ this restriction is not necessary because   $t^{-\frac{4-12\theta+n}{2(1-\theta)}} \rightarrow 0$  when $t  \rightarrow \infty$,  for any $\theta \in (0,\frac{1}{2})$.}
\end{rem}

Now we want to obtain an estimate for $ F_1(t,\cdot) $ on the region $ |\xi|\leq \eta^3 $. Initially, from \eqref{defieta1} we may see that 
\begin{align}\label{eqF4}
\int_{|\xi|\leq \eta^3} |F_1(t,\xi)|^2 d\xi &= \int_{|\xi|\leq \eta^3}   e^{-\frac{2|\xi|^2}{\log(1+|\xi|^{2\theta})}t} |R(\xi)|^2 |\hat{u}_1(\xi)|^2 d\xi \approx  \int_{|\xi| \leq \eta^3}  e^{-\log(1+|\xi|^{2-2\theta})t} |R(\xi)|^2 |\hat{u}_1(\xi)|^2 d\xi  \nonumber \\ 
&\leq \|u_1\|_{1}^2\int_{|\xi|\leq \eta^3}  e^{-\log(1+|\xi|^{2-2\theta})t} |R(\xi)|^2 d\xi
\end{align}
\noindent
Here, the function $R(r)$ is bounded on the low frequency zone for $0 < \theta \leq \frac{1}{3}$, because of  
\begin{equation}\label{ike-10}
\lim _{r \rightarrow +0} R(r) = \left\{\begin{matrix}
0 & \text{ for  } \,\, 0< \theta < \frac{1}{3},\\ 
 4 & \text{for }\,\, \theta =\frac{1}{3} .
\end{matrix}\right. 
\end{equation}
\noindent
Therefore, for $0 < \theta \leq \frac{1}{3} $ and $n \geq 1$,  from \eqref{eqF4}  and  \eqref{ike-10} we may conclude the existence of  a positive constant $C$  such that
\begin{align}\label{F_1theta<1/3}
\int_{|\xi|\leq \eta^3} |F_1(t,\xi)|^2 d\xi &\leq  C\|u_1\|_{1}^2\int_{|\xi|\leq \eta^3}  e^{-\log(1+|\xi|^{2-2\theta})t}  d\xi \nonumber \\
&= C\|u_1\|_{1}^2 \omega_n \int_{0}^{\eta^3}  (1+r^{2-2\theta })^{-t } r^{n-1}  dr \nonumber \\
&\sim   \|u_1\|_{1}^2 t^{-\frac{n}{2(1-\theta)}}, \quad t\gg 1. 
\end{align}

However, in the case of $ 0\leq \theta \leq \frac{5}{12} $ we also notice that the function  $R(r)\sqrt{r}$ in the region $|\xi| \leq \eta^3$ is bounded, because 
$$ \lim_{r \rightarrow 0} \sqrt{r}R(r) =  \left\{\begin{matrix}
0 & \text{ for }  0<  \theta < \frac{5}{12}, \\ 
4 & \text{ for } \,\, \theta =  \frac{5}{12}.
\end{matrix}\right. $$  
In particular, $R(r)\sqrt{r}$ in the region $r=|\xi| \leq \eta^3$ is bounded for $\frac{1}{4} < \theta \leq \frac{5}{12}$. Thus, from \eqref{eqF4}, in the case of $ \frac{1}{4} < \theta \leq \frac{5}{12} $ and $n \geq 2 $,  one can obtain
\begin{align}\label{F_1ngrande}
\int_{|\xi|\leq \eta^3}  |F_1(t,\xi)|^2 d\xi &\leq \omega_n \|u_1\|_{1}^2  \int_{0}^{\eta^3} (1+r^{2-2\theta})^{-t}|R(r)|^2 r^{n-1} dr \nonumber \\
&\leq C  \|u_1\|_{1}^2 \int_{0}^{\eta^3} (1+r^{2-2\theta})^{-t} r^{n-2} dr \nonumber \\
&\leq C \|u_1\|_{1}^2 t^{-\frac{n-1}{2(1-\theta)}}, \quad t\gg 1. 
\end{align}
\noindent
Similarly to the way used to obtain estimates for $F_1(t,\cdot) $ one can arrive at the following estimates for $F_2(t,\cdot)$:
\begin{align}\label{F_2decay}
\int_{|\xi|\leq \eta^3}  |F_2(t,\xi)|^2 d\xi \leq \left\{\begin{matrix}
 \frac{C}{\theta} \|u_1\|_1^2 t^{-\frac{n}{2\theta}}  &\text{ for  } n\geq 1 \text{ and } 0< \theta \leq \frac{1}{3}, \quad t\gg 1,\\ 
 \frac{C}{\theta} \|u_1\|_1^2 t^{-\frac{n-1}{2\theta}}  &\text{ for  } n\geq 2 \text{ and } \frac{1}{4}<  \theta \leq \frac{5}{12},  \quad t\gg 1.
\end{matrix}\right.
\end{align}


Let us estimate the $L^2$-norm of $ F_6(t,\xi)$ at the final stage in this subsection 3.1.1. To do that we need to analyze the function $ g(s)$ given by \eqref{defig}. Note that it is easy to see that 
\begin{equation}\label{glimitada}
1 \leq g(s)\leq 2
\end{equation}
for $s \leq \delta$, and its derivative is given by 
 \begin{align*}
  g'(s) = \frac{1}{2\sqrt{1- \frac{4s^6}{\log^2 (1+s^{6\theta})}} } \left(  \frac{ 48 \theta s^{6\theta+5}}{(1+s^{6\theta}) \log^3(1+s^{6\theta })} - \frac{24s^5}{\log^2(1+s^{6\theta})} \right).
 \end{align*}
Then, for $ \theta \in [0,\frac{5}{12} ] $, the function $ g'(s)$ is bounded on the interval  $0< s \leq \eta  $. In fact, the limits 
\begin{align*}
\lim _{s\rightarrow +0}\frac{  s^{6\theta+5}}{(1+s^{6\theta}) \log^3(1+s^{6\theta })}  
\qquad 
\mbox{and } \qquad 
\lim _{s \rightarrow +0}\frac{s^5}{\log^2(1+s^{6\theta})}
\end{align*}
are finite because of $ 0 \leq \theta \leq \frac{5}{12} $. It should be mentioned 
 that the same does note happen on the zone  $\eta<   s < \delta $ because
$$\lim_{s \to \delta-0}\left ( \sqrt{1- \frac{4s^6}{\log^2 (1+s^{6\theta})}} \right )^{-1} = +\infty$$
(See \eqref{definicao1delta}--\eqref{definicao2delta}\,). Recall again that for $\theta \in (0,1/2)$
\[\lim_{s \to +0}\frac{s^6}{\log^2(1+s^{6\theta})} = 0.\]  
\noindent
By summarizing above facts,  there exists a constant $K>0$ depending on $\theta \in [0,\frac{12}{5}]$ and $\eta > 0$ such that for all $ s \in [0,\eta]$ it holds that
\begin{align*}
 |g'(s)| \leq K .
\end{align*}
In particular, for $ |\xi| \in [0, \eta ^3] $, we have  $ \sqrt[3]{|\xi|} \in [0,\eta] $ and $ \alpha (t,\xi) \sqrt[3]{|\xi|} \in [0,\eta] $. Thus 
\begin{align}\label{infoderivadadeg1}
 |g'(\alpha \sqrt[3]{|\xi|})| \leq K , \quad |\xi| \leq \eta ^3. 
\end{align}
From \eqref{glimitada} and \eqref{infoderivadadeg1}, for $ 0< \theta\leq \frac{5}{12} $ and $n \geq 1$ we can estimate the $L^2$-norm of $F_6(t,\cdot)$ as follows: 
\begin{align}\label{AF_6}
\int _{|\xi| \leq \eta ^3} |F_6(t,\xi) |^2 d\xi &=  \frac{1}{4} t^{2} \int _{|\xi| \leq \eta ^3}  e^{ - 2\frac{t\log(1+|\xi|^{2\theta})}{2} g\left ( \alpha \sqrt[3]{|\xi|} \right ) } \frac{\log^2(1+|\xi|^{2\theta}) |\xi|^{\frac{2}{3}}  }{ (\lambda_+-\lambda _-)^2} | g'\left ( \alpha \sqrt[3]{|\xi|} \right )|^2 |\hat{u}_1(\xi)|^2 d\xi \nonumber \\
&\leq C t^{2}  \|u_1\|_1^2 \int _{|\xi| \leq \eta ^3}  e^{ - t\log(1+|\xi|^{2\theta}) }  |\xi|^{\frac{2}{3}} d\xi \nonumber \\
&= C \omega _n t^{2}  \|u_1\|_1^2 \int _{0}^{\eta ^3}   (1+r^{2\theta})^{-t}   r^{n-\frac{1}{3}}  dr \nonumber \\
&\sim  \frac{1}{\theta} t^{2}  \|u_1\|_1^2 t^{-\frac{n+\frac{2}{3}}{2\theta}}\nonumber \\
&= \frac{1}{\theta}   \|u_1\|_1^2 t^{-\frac{n-4\theta +\frac{2}{3}}{2\theta}}, \quad t\gg1.
\end{align}

As a result one can conclude the following Propositions. In that case, it is essential whether the factor $1/\theta$ can be included or not in the final estimates as the coefficient. 
\begin{pro}\label{Prop3.1}  Let $n=1$, $0 < \theta \leq \frac{1}{3}$, and $\varphi(t,\xi)$ be given by \eqref{defivarphi1}. If $u_1 \in L^{1,2\theta}({\bf R})$,  then 
\begin{align*}
\int _{|\xi|\leq \eta^{3}} \hspace{-0.2cm} |\hat{u}(t,\xi)-\varphi(t,\xi)|^2d\xi \leq \left\{\begin{matrix}
C(\|u_1\|_{1}^2+\|u_1\|_{L^{1,2\theta}}^2) \left ( t^{-\frac{1}{2(1-\theta)}} + \displaystyle{\frac{1}{\theta}}t^{-\frac{1}{2\theta}}\right ), & \text{ if } 0< \theta \leq \frac{1}{6}, \\ 
 C(\|u_1\|_{1}^2+\|u_1\|_{L^{1,2\theta}}^2) \left ( t^{-\frac{1}{2(1-\theta)}} + \displaystyle{\frac{1}{\theta}}t^{-\frac{\frac{5}{3}-4\theta}{2\theta}} \right ), & \text{ if } \frac{1}{6} < \theta \leq \frac{1}{3}
\end{matrix}\right. ,
\end{align*}
for $t\gg 1 $. 
\end{pro} 
{\it{Proof.}} The proof is obtained by choosing the slowest estimates as $t \to \infty$ among \eqref{estimativaF_3}, \eqref{estimativaF_5}, \eqref{F4rate}, \eqref{F_1theta<1/3}, \eqref{F_2decay} and \eqref{AF_6}. Note that the case $1/6 <\theta \leq 1/3$ is coming from the relation such that $ \frac{\frac{5}{3}-4\theta}{2\theta}\leq \frac{n}{2\theta} $ with $n = 1$. 
\hfill
$\Box$

\begin{pro}\label{Prop3.2} Let $n\geq 2$, $ 0 < \theta \leq \frac{5}{12} $ and $\varphi(t,\xi)$ be given by \eqref{defivarphi1}.  If $  u_1 \in L^{1,2\theta}({\bf R}^n) $,  then 
\begin{align*}
\int _{|\xi|\leq \eta^{3}} \hspace{-0.2cm} |\hat{u}(t,\xi)-\varphi(t,\xi)|^2d\xi \leq \left\{\begin{matrix}
C(\|u_1\|_{1}^2+\|u_1\|_{L^{1,2\theta}}^2) \left ( t^{-\frac{n}{2(1-\theta)}} + \displaystyle{\frac{1}{\theta}} t^{-\frac{n}{2\theta}}\right ), & \text{ if } 0< \theta \leq \frac{1}{6}, \\ 
 C(\|u_1\|_{1}^2+\|u_1\|_{L^{1,2\theta}}^2) \left ( t^{-\frac{n}{2(1-\theta)}} + \displaystyle{\frac{1}{\theta}} t^{-\frac{n-4\theta +\frac{2}{3}}{2\theta}} \right ), & \text{ if } \frac{1}{6} < \theta \leq \frac{1}{3},\\ 
C(\|u_1\|_{1}^2+\|u_1\|_{L^{1,2\theta}}^2) \left ( t^{-\frac{n-1}{2(1-\theta)}} + \displaystyle{\frac{1}{\theta}} t^{-\frac{n-1}{2\theta}} \right ),  & \text{ if } \frac{1}{3} < \theta \leq \frac{5}{12}
\end{matrix}\right.
\end{align*}
for $t \gg 1$.
\end{pro}
{\it{Proof.}} We may conclude this result by comparing the estimates \eqref{estimativaF_3}, \eqref{estimativaF_5}, \eqref{F4rate}, \eqref{F_1theta<1/3}, \eqref{F_1ngrande}, \eqref{F_2decay} and \eqref{AF_6}. Note that the case $1/6 <\theta \leq 1/3$ is coming from the relation such that $ \frac{n-4\theta+\frac{2}{3}}{2\theta}\leq \frac{n}{2\theta} $ with $n \geq 2$.
\hfill
$\Box$

\subsubsection{Estimates on the middle-frequency zone	$\eta^3 \leq |\xi| \leq \delta$}
We call the zone  $\eta^3 \leq  |\xi|\leq \delta $ the middle frequency  because the characteristics roots given by \eqref{defilambda} are real on this zone, and therefore the solution of \eqref{eqnfourier}-\eqref{initialfourier} is given by
\begin{equation*}
\hat{u}(t,\xi) = e^{-t\frac{\log(1+|\xi|^{2\theta})}{2}} \frac{\sinh (C(\xi)t)}{2C(\xi)} \hat{u}_1(\xi),
\end{equation*}
where
\begin{equation*}
C(\xi)=\frac{ \sqrt{\log^2(1+|\xi|^{2\theta})-4|\xi|^2} }{2}.
\end{equation*}

We remember that $ \eta $ is defined in \eqref{etainfimo}. Since the function $ \vert\xi\vert \mapsto  \frac{|\xi|^{4-4\theta}}{|\xi|^2} $ is increasing for $ 0< \theta < \frac{1}{2} $, we may observe that
 \begin{equation}\label{eta3supremo}
 \eta ^3 = \sup \{ \alpha >0; \frac{|\xi|^{4-4\theta}}{|\xi|^2} \leq \frac{1}{25^3} \text{ for } 0<| \xi|\leq \alpha  \} .
 \end{equation}

\begin{lem} \label{lemmiddle1}
There exists $ \beta =\beta(\theta)$, $ 0< \beta \leq \eta^3 $, such that 
\begin{align*}
\frac{2}{25^3} \log^2(1+|\xi|^{2\theta}) \geq 4 |\xi|^2 \text{ for } |\xi|\leq \beta \\
\frac{2}{25^3} \log^2(1+|\xi|^{2\theta}) \leq 4 |\xi|^2 \text{ for } |\xi|\geq \beta. 
\end{align*}
\end{lem}
{\it Proof.}\
The argument used to prove the existence of $\delta$ as in Lemma \ref{lemadelta} can be also used to prove the existence of $ \beta =\beta(\theta) \in (0,1) $, which satisfies both conclusions of this lemma. 
So, it suffices to check that $ \beta \leq \eta^3. $

From Remark \ref{rem3.1}, we know that $ \log^2(1+|\xi|^{2\theta}) \leq |\xi|^{4\theta} $, for $ |\xi|\leq 1 $. Thus, if $|\xi|\leq \beta$, we have 
$$  \frac{2}{25^3} |\xi|^{4\theta} \geq \frac{2}{25^3} \log^2(1+|\xi|^{2\theta}) \geq 4 |\xi|^2 .$$
This implies 
$$ 2 \times 25^3 |\xi|^{4-4\theta}  \leq |\xi|^2, \qquad  |\xi|  \leq \beta.$$
and the condition  $ \frac{|\xi|^{4-4\theta}}{|\xi|^2} \leq \frac{1}{25^3} $ is satisfied for $|\xi|\leq \beta$. Therefore, one has $\beta \leq \eta^3 $ from \eqref{eta3supremo}. 
\hfill
$\Box$

In other words,  Lemma \ref{lemmiddle1} tells us that 
 $$ \log^2 (1+|\xi|^{2\theta}) -4 |\xi|^2 < \frac{25^3-2}{25^3} \log^2(1+|\xi|^{2\theta})  \quad \mbox{ for}  \quad   |\xi| \geq \beta  $$
 and in particular, the definition of $C(\xi)$ implies that 
$$ 0 <  2 C(\xi)=\sqrt{\log^2 (1+|\xi|^{2\theta}) -4 |\xi|^2} < \frac{\sqrt{25^3-2}}{\sqrt{25^3}} \log (1+|\xi|^{2\theta})  \quad \mbox{ for}  \quad  \beta \leq |\xi| <\delta  . $$
Therefore, if $ \eta ^3 \leq |\xi| < \delta $ one has 
\begin{equation}\label{eqbeta}
- \log (1+|\xi|^{2\theta})  + 2 C(\xi)< \Big(\frac{\sqrt{25^3-2}}{\sqrt{25^3}}-1\Big) \log (1+|\xi|^{2\theta}) = - c \log (1+|\xi|^{2\theta})
\end{equation}
with  $0<c<1 $ a constant, due to the fact that $ \beta \leq \eta^3 $. 
\vspace{0.2cm}

Now, from Lemma \ref{lemmahiperbolicsine} and inequality \eqref{eqbeta} we can prove the exponential decay for the $L^2$-norm of  $\hat{u}(t,\cdot)$ on the middle frequency zone as follows:
\begin{align}\label{eq3.47}
\int_{\eta^3 \leq|\xi| < \delta} |\hat{u}(t,\xi)|^2d\xi &= \int_{\eta^3 \leq|\xi| < \delta} e^{-t \log(1+|\xi|^{2\theta})} \frac{\sinh^2 (C(\xi)t)}{4(C(\xi))^2} |\hat{u}_1(\xi)|^2d\xi \nonumber \\
&\leq  \frac{K^{2}}{4}t^2 \int_{\eta ^3\leq|\xi|  \leq \delta} e^{-t\log(1+|\xi|^{2\theta}) +2C(\xi)t}  |\hat{u}_1(\xi)|^2d\xi \nonumber \\
&\leq  K^{2}t^2 \int_{\eta^3\leq|\xi|\leq \delta} e^{-ct\log(1+|\xi|^{2\theta})}   |\hat{u}_1(\xi)|^2d\xi \quad ( c>0)\nonumber \\
&= K^{2}t^2 \int_{\eta^3\leq|\xi|\leq \delta} (1+|\xi|^{2\theta})^{-ct}  |\hat{u}_1(\xi)|^2d\xi \nonumber \\
&\leq K^{2} \omega_nt^2 \|u_1\|_1^2 \int_{\eta^3}^{\delta} (1+r^{2\theta})^{-ct} r^{n-1}dr \nonumber \\
&\leq C t^2 (1+\eta^{6\theta})^{-ct}\|u_1\|_1^2, \qquad t \gg 1,
\end{align}
with $C$ a positive constant depending on the space dimension $n$ and $c>0$ a constant given in  \eqref{eqbeta}.

\subsection{Estimates on the high-frequency zone $|\xi| \geq \delta$}

On the high frequency zone $|\xi| > \delta$  the characteristics roots are complex and the solution of \eqref{eqnfourier}-\eqref{initialfourier} is given by
\begin{equation*}
\hat{u}(t,\xi)= \frac{1}{b(\xi)} e^{-a(\xi)t}  \sin(b(\xi)t)\hat{u}_1(\xi)
\end{equation*}
where
\begin{equation*}
a(\xi) = \frac{ \log (1+|\xi|^{2\theta})}{2}, \quad b(\xi)= \frac{ \sqrt{ 4 |\xi|^2 -\log ^2(1+|\xi|^{2\theta})  }}{2}.
\end{equation*}
We know that $ \vert\sin a\vert \leq a $ for all $a\geq 0 $. Then $\vert \displaystyle{\frac{\sin(b(\xi)t)}{b(\xi)}}\vert \leq t$ for all $t\geq 0$, and so one has 
\begin{align}\label{eq3.49}
\int _{|\xi|  > \delta} |\hat{u}(t,\xi)|^2d\xi &=  \int _{|\xi| > \delta} (1+|\xi|^{2\theta})^{-t} \frac{\sin^2(b(\xi)t)}{b(\xi)^2} |\hat{u}_1(\xi)|^2 d\xi \nonumber \\
&\leq t^2 \|u_1\|_1^2 \int _{|\xi|\geq\delta} (1+|\xi|^{2\theta})^{-t}   d\xi \nonumber \\
&=\omega _n t^2 \|u_1\|_1^2 \int _{\delta}^{1} (1+r^{2\theta})^{-t} r^{n-1}   dr  + \omega _n t^2 \|u_1\|_1^2 \int _{1}^{\infty} (1+r^{2\theta})^{-t} r^{n-1}   dr \nonumber \\
&\sim \|u_1\|_1^2 t^2 \Big((1+\delta^{2\theta})^{-t} +  \frac{2^{-t}}{t-1} \Big), \qquad t \gg 1. 
\end{align}
The last inequality is obtained by using Lemma \ref{infit}.


\subsection{Proof of Theorem \ref{maintheorem}}

Now by combining Propositions \ref{Prop3.1}, \ref{Prop3.2}, \eqref{eq3.47}, and \eqref{eq3.49} one can prove our main Theorem \ref{maintheorem}.\\

{\it Proof of Theorem \ref{maintheorem}.}\
We first note that $ \log(1+|\xi|^{2\theta}) \leq |\xi|^{2\theta} $ for all $ \xi \in {\bf R}^n$,  which implies  $ \displaystyle{\frac{|\xi|^{2\theta}}{\log (1+|\xi|^{2\theta})}} \geq 1 $. Then, one can get the next  estimate for $ t\gg 1 $ on the zone of high frequency $|\xi| \geq \eta^3$ as follows:

\begin{align}\label{eq3.50}
\int _{|\xi|\geq \eta^3} |\varphi(t,\xi)|^2d\xi &\leq P_1^2 \int _{|\xi|\geq \eta^3}  \frac{e^{-\frac{2|\xi|^2}{\log(1+|\xi|^{2\theta})}t}}{\log ^2(1+|\xi|^{2\theta})} d\xi + P_1^2 \int _{|\xi|\geq \eta^3}  \frac{e^{-\log(1+|\xi|^{2\theta})t}}{\log^2(1+|\xi|^{2\theta})} d\xi \nonumber \\
&=P_1^2 \int _{|\xi|\geq \eta^3}  \frac{e^{-\frac{2|\xi|^{2\theta}  |\xi|^{2-2\theta}}{\log(1+|\xi|^{2\theta})}t}}{\log ^2(1+|\xi|^{2\theta})} d\xi + P_1^2 \int _{|\xi|\geq \eta^3}  \frac{(1+|\xi|^{2\theta})^{-t}}{\log^2(1+|\xi|^{2\theta})} d\xi \nonumber \\
&\leq  P_1^2 \int _{|\xi|\geq \eta^3}  \frac{e^{-2t|\xi|^{2-2\theta}} }{\log ^2(1+|\xi|^{2\theta})} d\xi + P_1^2 \int _{|\xi|\geq \eta^3}  \frac{(1+|\xi|^{2\theta})^{-t}}{\log^2(1+|\xi|^{2\theta})} d\xi \nonumber \\
&\leq \frac{P_1^2}{\log ^2(1+\eta^{6\theta})} \int_{|\xi|\geq \eta^3} e^{-2t|\xi|^{2-2\theta}} d\xi + \frac{P_1^2 \omega_n}{\log ^2(1+\eta^{6\theta})} \int_{\eta^3}^{1} (1+r^{2\theta})^{-t} r^{n-1} dr \nonumber \\
&+ \frac{P_1^2 \omega_n}{(\log 2)^{2}} \int_{1}^{\infty} (1+r^{2\theta})^{-t} r^{n-1} dr \nonumber \\
&\leq \frac{P_1^2}{\log ^2(1+\eta^{6\theta})} e^{-t\eta^{6-6\theta}}\int_{|\xi|\geq \eta^3} e^{-|\xi|^{2-2\theta}} d\xi + \frac{P_1^2 \omega_n}{\log ^2(1+\eta^{6\theta})} \int_{\eta^3}^{1} (1+r^{2\theta})^{-t} r^{n-1} dr \nonumber \\
&+ \frac{P_1^2 \omega_n}{(\log 2)^{2}} \int_{1}^{\infty} (1+r^{2\theta})^{-t} r^{n-1} dr \nonumber \\
&\leq C P_1^2 \Big( e^{-t\eta^{6-6\theta}} +  (1+\eta^{6\theta})^{-t} + \frac{2^{-t}}{t-1} \Big),\quad  t\gg 1. 
\end{align}

Now, it follows from the Plancherel Theorem that 
$$ \int _{{\bf R}^n}| u(t,x) - {\cal F}_{\xi\rightarrow x}^{-1}(\varphi (t,\xi))(x)|^2 dx = \int _{{\bf R}^n} |\hat{u}(t,\xi)-\varphi(t,\xi)|^2d\xi $$
for $t\geq 0$. Furthermore, one has
\begin{align}\label{eq3.51}
\int _{{\bf R}^n} |\hat{u}(t,\xi)-\varphi(t,\xi)|^2d\xi &\leq \int _{|\xi|\leq \eta^3} |\hat{u}(t,\xi)-\varphi(t,\xi)|^2d\xi + \int _{|\xi|\geq \eta^3} |\hat{u}(t,\xi)|^2d\xi + \int _{|\xi|\geq \eta^3} |\varphi(t,\xi)|^2d\xi 
\end{align}
for $t >  0$.

From \eqref{eq3.47} and \eqref{eq3.49}, we know that the $L^2$-estimates  on the zone $|\xi|\geq\eta^3 $ to $\hat{u}(t,\xi)$ are of exponential type.  The estimate to $\varphi(t,\xi) $ on $|\xi|\geq \eta^3 $ obtained in \eqref{eq3.50} is also faster than those obtained in Propositions \ref{Prop3.1} and \ref{Prop3.2}. The result of Theorem \ref{maintheorem} follows by combining  Propositions \ref{Prop3.1} and \ref{Prop3.2}, with inequalities \eqref{eq3.47}, \eqref{eq3.49}, \eqref{eq3.50} and \eqref{eq3.51}.
\hfill
$\Box$



\section{Optimality of the decay  rates}

Our goal in this section is to prove Theorem \ref{main-teo}, which shows the optimal decay rates in time depending on the dimension $n$ to the solution of the problem \eqref{eqn}-\eqref{initial}. From \eqref{defivarphi1}, we have
\begin{equation}\label{varphiaux}
\varphi(t,\xi) = \varphi _1(t,\xi) - \varphi_2(t,\xi), \quad t\geq 0, \quad \xi \in {\bf R}^n,
\end{equation}
 where
\begin{align*}
&\varphi_1(t,\xi):= \frac{e^{-\frac{|\xi|^2}{\log(1+|\xi|^{2\theta})}t}}{\log (1+|\xi|^{2\theta})} P_1, \quad \varphi_2(t,\xi) := \frac{e^{-\log(1+|\xi|^{2\theta})t}}{\log(1+|\xi|^{2\theta})}P_1. 
\end{align*}

\begin{lem}\label{lemma4.1} Let $ n=1 $ with $ 0< \theta < \frac{1}{4} $ and $ n\geq 2 $ with $ 0<  \theta \leq \frac{5}{12} $. If $ u_1 \in L^{1}({\bf R}^n) $, then 
$$ C_1 P_1^2 t^{-\frac{n-4\theta}{2(1-\theta)}}  \leq \int _{{\bf R}^n } |\varphi (t,\xi)|^2 d\xi \leq C_2 P_1^2 \left ( t^{-\frac{n-4\theta}{2(1-\theta)}} +  \frac{1}{\theta}  t^{- \frac{n-4\theta}{2\theta}} \right )
 , \quad t\gg 1,$$
where the constants $C_1, C_2$ depend only on $ \theta $ and $n$.  
\end{lem} 
{\it{Proof.}}

 First we note that
\begin{align}\label{estimatevarphi}
\int _{{\bf R}^n } |\varphi (t,\xi)|^2 d\xi &\leq  2 \int _{{\bf R}^n}   |\varphi_1(t,\xi)|^2  d\xi + 2 \int _{{\bf R}^n}  |\varphi_2(t,\xi)|^2 d\xi \nonumber \\
&= 2 \int _{|\xi| \leq \eta }   |\varphi_1(t,\xi)|^2  d\xi + 2 \int _{|\xi| \geq \eta }   |\varphi_1(t,\xi)|^2   d\xi \nonumber \\
&+2 \int _{|\xi| \leq \eta } |\varphi_2(t,\xi)|^2 d\xi + 2 \int _{|\xi| \geq \eta }  |\varphi_2(t,\xi)|^2 d\xi , \quad t > 0.
\end{align}
By using the equivalences obtained in Remark \ref{rem3.1}, we have
\begin{align}\label{1decayu}
 \int _{|\xi| \leq \eta}   |\varphi_1(t,\xi)|^2  d\xi &\approx P_1^2 \int _{|\xi| \leq \eta }   \frac{  e^{-t\log(1+|\xi|^{2-2\theta})}}{\log^2 (1+|\xi|^{2\theta})}  d\xi \leq P_1^2 \int _{|\xi| \leq \eta }   \frac{  (1+|\xi|^{2-2\theta})^{-t} }{\log^2 (1+|\xi|^{2\theta})}  d\xi \nonumber \\
&=\omega _n P_1^2 \int _{0}^{\eta}   \frac{(1+r^{2-2\theta})^{-t} }{\log^2 (1+r^{2\theta})} r^{n-1} dr =\omega _n  P_1^2 \int _{0}^{\eta}   \frac{(1+r^{2-2\theta})^{-t} }{\log^2 (1+r^{2\theta})} r^{n-1-4\theta} r^{4\theta} dr \nonumber \\
&\leq 4\omega _n P_1^2 \int _{0}^{\eta}   \frac{(1+r^{2-2\theta})^{-t} }{r^{4\theta}} r^{n-1-4\theta} r^{4\theta} dr = 4\omega _n P_1^2 \int _{0}^{\eta}  (1+r^{2-2\theta})^{-t}  r^{n-1-4\theta}  dr \nonumber \\
&\leq C P_1^2 t^{-\frac{n-4\theta}{2(1-\theta)}}, \quad t\gg 1. 
\end{align}
The last decay estimate is obtained from Lemma \ref{lema2-2theta} since $ n-4\theta >0$.
In the same way, by using Lemma \ref{lema2theta} for $n-4\theta>0$, we have the next estimate. 
\begin{align}\label{2decayu}
\int _{|\xi| \leq \eta } |\varphi_2(t,\xi)|^2 d\xi &=  P_1^2 \int _{|\xi| \leq \eta }  \frac{e^{-2t\log(1+|\xi|^{2\theta})} }{\log^2 (1+|\xi|^{2\theta})}d\xi \leq \omega _n P_1^2 \int_{0}^{\eta} \frac{(1+r^{2\theta})^{-t} }{\log ^2 (1+r^{2\theta})} r^{n-1} dr \nonumber \\
&\leq 4\omega _n P_1^2 \int_{0}^{\eta} (1+r^{2\theta})^{-t}  r^{n-1 -4\theta} dr \leq C \frac{1}{\theta} P_1^2 t^{- \frac{n-4\theta}{2\theta}}, \,\, t\gg 1. 
\end{align}
Further, from \eqref{eq3.50}, the $L^2$-estimate to $ \varphi (t,\xi) $ on the zone $|\xi| \geq \eta $  is of  exponential type, because $|\xi| \geq \eta $  implies that $|\xi| \geq \eta^3 $. 
Therefore, there exists a constant $C>0$ such that
\begin{align*}
\int _{{\bf R}^n } |\varphi (t,\xi)|^2 d\xi \leq C P_1^2 \left ( t^{-\frac{n-4\theta}{2(1-\theta)}} +  \frac{1}{\theta}  t^{- \frac{n-4\theta}{2\theta}} \right ) , \quad t\gg 1 ,  
\end{align*}
due to \eqref{estimatevarphi}, \eqref{1decayu} and \eqref{2decayu}.
 \vspace{0.2cm}
 
 In order to prove the estimate from bellow, from Remark \ref{rem3.1} we have
\begin{align}\label{lowestivarphi1}
\int_{{\bf R}^n} |\varphi _1(t,\xi)|^2 d\xi &\geq \int_{ |\xi|\leq \eta } |\varphi _1(t,\xi)|^2 d\xi \approx  P_1^2 \int _{|\xi| \leq \eta }   \frac{  e^{-t\log(1+|\xi|^{2-2\theta})}}{\log^2 (1+|\xi|^{2\theta})}  d\xi \nonumber \\
&= \omega _n  P_1^2 \int _{0 }^{\eta}   \frac{  e^{-t\log(1+r^{2-2\theta})}}{\log^2 (1+r^{2\theta})} r^{n-1}  dr \geq C\omega _n  P_1^2 \int _{0 }^{\eta}   \frac{  e^{-t\log(1+r^{2-2\theta})}}{r^{4\theta}} r^{n-1}  dr \nonumber \\
&= C\omega _n  P_1^2 \int _{0 }^{\eta} (1+r^{2-2\theta})^{-t}  r^{n-1-4\theta}  dr \nonumber \\
&\geq C P_1^2 t^{- \frac{n-4\theta}{2(1-\theta)}}, 
\end{align} 
because of $n-4\theta > 0$, due to Remark \ref{estimativaporbaixoIp}, where $C > 0$ is a generous constant.
\vspace{0.2cm}
 We also notice that $|\varphi_1(t,\xi)| \leq |\varphi(t,\xi)| + |\varphi_2(t,\xi)|$ and, from Young's inequality, $|\varphi_1(t,\xi)|^2 \leq 2|\varphi(t,\xi)|^2 +2 |\varphi_2(t,\xi)|^2 $. Thus, 
\begin{align*}
|\varphi(t,\xi)|^2 \geq \frac{1}{2}|\varphi_1(t,\xi)|^2 - |\varphi_2(t,\xi)|^2 , t\geq 0, \xi \in {\bf R}^n. 
\end{align*}
Then, from \eqref{lowestivarphi1} and \eqref{2decayu}, we have 
 \begin{align}\label{E1}
 \int _{|\xi| \leq \eta}|\varphi(t,\xi)|^2 d\xi &\geq \frac{1}{2} \int _{|\xi| \leq \eta} |\varphi_1(t,\xi)|^2 d\xi - \int _{|\xi| \leq \eta} |\varphi_2(t,\xi)|^2 d\xi \nonumber\\
&\geq K_1 P_1^2 t^{- \frac{n-4\theta}{2(1-\theta)}} - K_2 \frac{1}{\theta} P_1^2 t^{- \frac{n-4\theta}{2\theta}} \nonumber \\
&= P_1^2  t^{- \frac{n-4\theta}{2(1-\theta)}} \left ( K_1 - K_2 \frac{1}{\theta} t^{- \frac{8\theta^2 -2\theta n +n -4\theta}{2\theta(1-\theta)}} \right ).
 \end{align}
Since $ 0<  \theta < \frac{1}{2} $  and $ n-4\theta > 0 $, one can conclude that  $ 8\theta^2 -2\theta n +n -4\theta >0 $. 
Therefore, it follows from \eqref{E1} that 
\begin{align*}
\int _{{\bf R}^n}|\varphi(t,\xi)|^2 d\xi \geq \int _{|\xi| \leq \eta}|\varphi(t,\xi)|^2 d\xi &\geq \frac{K_1}{2} P_1^2  t^{- \frac{n-4\theta}{2(1-\theta)}} , \quad t \gg 1. 
\end{align*}
These arguments imply the desired estimate for $\varphi(t,\xi)$.
\hfill
$\Box$
\vspace{0.2cm}

The above arguments do not hold for $ n=1$ and $ \frac{1}{4} \leq \theta \leq \frac{1}{3} $, because the integrals  $$ \int_{0}^{\eta} \frac{(1+r^{2\theta})^{-t} }{\log ^2 (1+r^{2\theta})} r^{n-1} dr, \quad \int_{0}^{\eta} \frac{(1+r^{2 -2\theta})^{-t} }{\log ^2 (1+r^{2\theta})} r^{n-1} dr $$
are divergent for all $t > 0$. For this reason, we need to estimate the $L^2$-norm of the function $\varphi(t,\xi)$ itself:  
$$ \varphi(t,\xi) = \frac{e^{-\frac{|\xi|^2}{\log(1+|\xi|^{2\theta})}t}}{\log (1+|\xi|^{2\theta})} P_1 -   \frac{e^{-\log(1+|\xi|^{2\theta})t}}{\log(1+|\xi|^{2\theta})}P_1 .$$

\begin{lem}\label{lemma4.2}Let $n=1$ and $ \theta > \frac{1}{4} $. If $ u_1 \in L^{1}({\bf R}) $, there exist constants $C_1, C_2 >0$ such that 
$$ C_1 P_1^2 t^{ \frac{4\theta -1}{2\theta }}\leq \int _{{\bf R}} |\varphi (t,\xi)|^2 d\xi \leq C_2 \frac{1}{4\theta-1} P_1^2 t^{ \frac{4\theta -1}{2\theta }} , \quad t\gg 1.$$
\end{lem}
{\it Proof.}
We first note that from \eqref{eq3.50} the $L^2$-norm of $ \varphi(t,\xi) $ decays  exponentially on the high frequency region $ |\xi|\geq \eta >\eta^3 $. So, in this proof it suffices to consider the integral only in the low frequency zone $0 < \vert\xi\vert \leq \eta$.
\vspace{0.2cm }

Now, we notice that
\begin{equation*}
-\log(1+r^{2\theta}) = \frac{r^2-\log^2(1+r^{2\theta})}{\log(1+r^{2\theta})} - \frac{r^2}{\log(1+r^{2\theta})},
\end{equation*}
so that one has
\begin{align}\label{A5.2}
\log(1+\vert\xi\vert^{2\theta})\varphi(t,\xi) &= P_{1}\left(e^{-\frac{|\xi|^2}{\log(1+|\xi|^{2\theta})}t} - e^{-\log(1+|\xi|^{2\theta})t}\right) \nonumber\\
&= P_{1}\left(e^{-\frac{|\xi|^2}{\log(1+|\xi|^{2\theta})}t} -  e^{t\frac{|\xi|^2-\log^2(1+|\xi|^{2\theta})}{\log(1+|\xi|^{2\theta})} -t \frac{|\xi|^2}{\log(1+|\xi|^{2\theta})}} \right)\nonumber \\
&= P_{1} e^{-\frac{|\xi|^2}{\log(1+|\xi|^{2\theta})}t} \left ( 1 - e^{-t\frac{\log^2(1+|\xi|^{2\theta})-|\xi|^2}{\log(1+|\xi|^{2\theta})}} \right ).
\end{align}
\noindent
Due to the fact that for $ 0\leq r \leq 1 $ we have $ \frac{1}{2}r^{2\theta} \leq \log(1+r^{2\theta}) \leq r^{2\theta} $,  thus one has 
\begin{equation}\label{A5.5}
 r^{2\theta} \frac{(1-4r^{2-4\theta})}{4}\leq \frac{\log^2(1+r^{2\theta}) - r^2 }{ \log(1+r^{2\theta})} \leq 2r^{2\theta}(1-r^{2-4\theta}).
\end{equation}
\noindent
Moreover, since $ \theta < \frac{1}{2} $ we have $ 2-4\theta >0 $. Therefore, there exists $ \beta = \beta(\theta)>0 $, with $ \beta \leq  \eta $ such that 
$$ 1- 4r^{2-4\theta} \geq \frac{1}{2} ,$$ 
for $ 0\leq r\leq \beta $. Thus, 
\begin{equation}\label{A5.6}
\frac{1}{2} \leq 1- 4r^{2-4\theta} \leq 1- r^{2-4\theta} \leq 1 
\end{equation}
for $ 0\leq r \leq \beta$. From \eqref{A5.5} and \eqref{A5.6} one can get
\begin{equation*}
\frac{1}{8} r^{2\theta} \leq \frac{\log^2(1+r^{2\theta}) - r^2 }{ \log(1+r^{2\theta})} \leq 2r^{2\theta},
\end{equation*}
for $ 0< r\leq \beta $.
This implies 
\begin{equation*}
1-e^{-\frac{1}{8} tr^{2\theta}}\leq 1-e^{-t \frac{\log^2(1+r^{2\theta}) - r^2 }{ \log(1+r^{2\theta})}} \leq 1-e^{-2tr^{2\theta}}
\end{equation*}
and
\begin{equation}\label{A5.9}
\frac{1-e^{-\frac{1}{8} tr^{2\theta}} }{r^{2\theta}}\leq \frac{1-e^{-t \frac{\log^2(1+r^{2\theta}) - r^2 }{ \log(1+r^{2\theta})}}}{\log(1+r^{2\theta})} \leq 2\frac{1-e^{-2tr^{2\theta}}}{r^{2\theta}}, 
\end{equation}
for $0< r\leq \beta $. Since 
$$ \lim _{\sigma \rightarrow 0} \frac{1-e^{-\sigma}}{\sigma} = 1,$$
there exists $ \alpha >0 $ such that $ \alpha  \leq  \beta  $ and 
\begin{equation}\label{A5.10}
\frac{1}{2}  \leq \frac{1-e^{-\sigma}}{\sigma} \leq \frac{3}{2},
\end{equation} 
for all $0< \sigma \leq \alpha $.
\vspace{0.2cm}
Based on these preparations let us prove the desired estimate for $\varphi(t,\xi)$.\\
\noindent
{\rm (1) The lower estimate of Lemma:} \\

For $ 0< r\leq \left ( \frac{8\alpha}{t} \right )^{\frac{1}{2\theta}} $ it holds that  $0< \sigma=\frac{1}{8} tr^{2\theta} \leq \alpha $. Applying estimate \eqref{A5.10} we get 
\begin{equation}\label{A5.11}
\frac{1}{2}  \leq \frac{1-e^{-\frac{1}{8} tr^{2\theta}} }{\frac{1}{8}t r^{2\theta}} \leq \frac{3}{2} .
\end{equation}
From \eqref{A5.9} and \eqref{A5.11}, for $0< r \leq    \left ( \frac{8\alpha}{t} \right )^{\frac{1}{2\theta}} $, it holds that
\begin{equation}\label{A5.12}
\frac{1-e^{-t \frac{\log^2(1+r^{2\theta}) - r^2 }{ \log(1+r^{2\theta})}}}{\log(1+r^{2\theta})}  \geq \frac{t}{16} .
\end{equation}
Let $ t_0> 0 $ be such that  $\left ( \frac{8\alpha}{t_0} \right )^{\frac{1}{2\theta}} \leq \alpha $, and consider  $ t \geq t_0 $. By combining \eqref{varphiaux} with \eqref{A5.2} and \eqref{A5.12},  since $\alpha \leq \beta \leq \eta $, we obtain 
\begin{align}\label{A5.13}
\int _{|\xi| \leq \eta} |\varphi (t,\xi)|^2 d\xi &= P_1^2 \int _{|\xi| \leq \eta} \left (  \frac{e^{-\frac{|\xi|^2}{\log(1+|\xi|^{2\theta})}t} - e^{-\log(1+|\xi|^{2\theta})t}  }{\log (1+|\xi|^{2\theta})} \right )^2 d\xi \nonumber \\ 
&= \omega _1 P_1^2 \int _{0}^{ \eta} e^{-\frac{2r^2}{\log(1+r^{2\theta})}t} \left ( \frac{1 - e^{-t\frac{\log^2(1+r^{2\theta})-r^2}{\log(1+r^{2\theta})}}}{\log (1+r^{2\theta})} \right )^2 dr \nonumber \\
&\geq \frac{\omega _1}{16^2}  P_1^2 t^2 \int _{0}^{ \left ( \frac{8\alpha}{t} \right )^{\frac{1}{2\theta}} } e^{-\frac{2r^2}{\log(1+r^{2\theta})}t} dr, \quad t\geq t_0. 
\end{align}
We also notice that 
\begin{equation}\label{A5.14}
2 \log (1+r^{2-2\theta})  \leq 2r^{2-2\theta} \leq  \frac{2r^2}{\log(1+r^{2\theta})} \leq 4 r^{2-2\theta} \leq 8 \log (1+r^{2-2\theta}) \quad 0<r\leq 1.
\end{equation}
Thus, from \eqref{A5.13} and \eqref{A5.14} one has 
\begin{align}\label{low-1}
\int _{|\xi| \leq \eta} |\varphi(t,\xi)|^2 d\xi &\geq \frac{\omega _1}{16^2}  P_1^2 t^2 \int _{0}^{ \left ( \frac{8\alpha}{t} \right )^{\frac{1}{2\theta}} } e^{-8t \log(1+r^{2-2\theta})} dr \nonumber \\
 &=  \frac{\omega _1}{16^2}  P_1^2 t^2 \int _{0}^{ \left ( \frac{8\alpha}{t} \right )^{\frac{1}{2\theta}} } (1+r^{2-2\theta})^{-8t} dr \nonumber \\
 &\geq \frac{\omega _1}{16^2}  P_1^2 t^2  \left ( 1+ \left ( \frac{8\alpha}{t} \right )^{\frac{2-2\theta}{2\theta}} \right )^{-8t}   \int _{0}^{ \left ( \frac{8\alpha}{t} \right )^{\frac{1}{2\theta}} } dr, \quad t\geq t_0.
\end{align}
Now we observe that  $1< \frac{2-2\theta}{2\theta} < 3 $ for  $  \frac{1}{4}< \theta < \frac{1}{2} $.  Then there exists $T\geq t_0$ such that 
\begin{equation} \label{low-2}
\frac{1}{2} \leq \left ( 1+ \left ( \frac{8\alpha}{t} \right )^{\frac{2-2\theta}{2\theta}} \right )^{-t} \leq \frac{3}{2}
\end{equation}
for all $t \geq T$, because of the fact that 
 \begin{equation*}
 \lim _{t \rightarrow +\infty}\left ( 1+ \frac{1}{t^q} \right )^{-t} =1 
 \end{equation*}
provided that $q > 1$. By combining  estimates \eqref{low-1} and \eqref{low-2} one can arrive at the desired 
 estimate from below such that
\begin{align*}
\int _{|\xi| \leq \eta} |\varphi(t,\xi)|^2 d\xi &\geq \frac{\omega _n}{16^2}  P_1^2 t^2  \left ( 1+ \left ( \frac{8\alpha}{t} \right )^{\frac{2-2\theta}{2\theta}} \right )^{-t}   \int _{0}^{ \left ( \frac{8\alpha}{t} \right )^{\frac{1}{2\theta}} } dr  \nonumber \\
&\geq \frac{\omega _n}{2 \times 16^2}  P_1^2 t^2  \int _{0}^{ \left ( \frac{8\alpha}{t} \right )^{\frac{1}{2\theta}} } dr  \nonumber \\
&= \frac{\omega _n}{2 \times 16^2}  P_1 t^2 \left ( \frac{8\alpha}{t} \right )^{\frac{1}{2\theta}} \nonumber \\
&= C P_1^2 t^{2}t^{-\frac{1}{2\theta}} \nonumber \\
&= C P_1^2  t^{\frac{4\theta-1}{2\theta}}, \quad t\geq T
\end{align*} 
with some constant $C = C_{\theta} > 0$.\\

\noindent 
{\rm (2) The upper estimate of Lemma:} \\
\noindent
From \eqref{A5.2}, we have 
\begin{align*}
\int _{|\xi| \leq \eta} |\varphi(t,\xi) |^2 d\xi &= \omega _1 P_1^2 \int _{0}^{ \eta} e^{-\frac{2r^2}{\log(1+r^{2\theta})}t} \left ( \frac{1 - e^{-t\frac{\log^2(1+r^{2\theta})-r^2}{\log(1+r^{2\theta})}}}{\log (1+r^{2\theta})} \right )^2 dr = A_{1}(t, \theta) + A_{2}(t, \theta),
 \end{align*}
where 
\begin{align}
A_1(t, \theta) &:= \omega _1 P_1^2 \int _{0}^{\left ( \frac{\alpha}{2t} \right )^{\frac{1}{2\theta}} } e^{-\frac{2r^2}{\log(1+r^{2\theta})}t} \left ( \frac{1 - e^{-t\frac{\log^2(1+r^{2\theta})-r^2}{\log(1+r^{2\theta})}}}{\log (1+r^{2\theta})} \right )^2 dr \label{defiI_1}, \\
A_2(t,  \theta) &:=  \omega _1 P_1^2 \int _{\left ( \frac{\alpha}{2t} \right )^{\frac{1}{2\theta}} }^{ \eta} e^{-\frac{2r^2}{\log(1+r^{2\theta})}t} \left ( \frac{1 - e^{-t\frac{\log^2(1+r^{2\theta})-r^2}{\log(1+r^{2\theta})}}}{\log (1+r^{2\theta})} \right )^2 dr, \label{defiI_2}
\end{align} 
 which holds for $t\geq t_0$.
 
Now, for  $ 0 < r \leq \left ( \frac{\alpha}{2t} \right )^{\frac{1}{2\theta}} $ by using  inequality \eqref{A5.10} we  have that 
\begin{align}\label{A5.19}
\frac{1-e^{-2tr^{2\theta}}}{r^{2\theta}}\leq 3t .
\end{align}
Thus for $0< r \leq    \left ( \frac{\alpha}{2t} \right )^{\frac{1}{2\theta}} $, by combining \eqref{A5.9} with \eqref{A5.19} it holds that
\begin{equation}\label{A5.20}
\frac{1-e^{-t \frac{\log^2(1+r^{2\theta}) - r^2 }{ \log(1+r^{2\theta})}}}{\log(1+r^{2\theta})}  \leq 6t , \quad t\geq t_0.
\end{equation}
The definition of $A_1(t,\theta)$ and the inequality \eqref{A5.20} imply that
\begin{align}\label{estimateA_1}
A_1(t, \theta) &=  \omega _1 P_1^2 \int _{0}^{\left ( \frac{\alpha}{2t} \right )^{\frac{1}{2\theta}} } e^{-\frac{2r^2}{\log(1+r^{2\theta})}t} \left ( \frac{1 - e^{-t\frac{\log^2(1+r^{2\theta})-r^2}{\log(1+r^{2\theta})}}}{\log (1+r^{2\theta})} \right )^2 dr \nonumber \\
&\leq  36 t^2 \omega _1 P_1^2 \int _{0}^{ \left ( \frac{\alpha}{2t} \right )^{\frac{1}{2\theta}} } e^{-\frac{2r^2}{\log(1+r^{2\theta})}t} dr  \leq  36 t^2 \omega _1 P_1^2 \int _{0}^{ \left ( \frac{\alpha}{2t} \right )^{\frac{1}{2\theta}} }  dr \nonumber \\ 
&= 36 t^2 \omega _1 P_1^2 \left ( \frac{\alpha}{2t} \right )^{\frac{1}{2\theta}}  = C P_1^2 t^2 t ^{- \frac{1}{2\theta}} \nonumber \\
&= C P_1^2 t^{\frac{4\theta-1}{2\theta}}, \quad t\gg 1 
\end{align}
with some generous constant $C > 0$. Note that the estimate given by \eqref{estimateA_1} is also holds for $ \theta = \frac{1}{4} $. \\

In order to estimate $ A_2(t,\theta)$, we use \eqref{A5.9} for $\theta >1/4$ to get the following estimate 
\begin{align}
A_2(t, \theta)&=  \omega _1 P_1^2 \int _{\left ( \frac{\alpha}{2t} \right )^{\frac{1}{2\theta}} }^{ \eta} e^{-\frac{2r^2}{\log(1+r^{2\theta})}t} \left ( \frac{1 - e^{-t\frac{\log^2(1+r^{2\theta})-r^2}{\log(1+r^{2\theta})}}}{\log (1+r^{2\theta})} \right )^2 dr \nonumber \\
 &\leq 4\omega _1 P_1^2 \int _{\left ( \frac{\alpha}{2t} \right )^{\frac{1}{2\theta}} }^{ \eta} e^{-\frac{2r^2}{\log(1+r^{2\theta})}t} \left ( \frac{1 - e^{-2t r^{2\theta}}}{r^{2\theta}} \right )^2 dr \nonumber \\ 
 &\leq 4\omega _1 P_1^2 \int _{\left ( \frac{\alpha}{2t} \right )^{\frac{1}{2\theta}} }^{ \eta} \frac{1}{r^{4\theta}} dr \nonumber \\
&= 4\omega _1 P_1^2  \frac{1}{1-4\theta} \left ( \eta ^{1-4\theta} - \left ( \frac{\alpha}{2t} \right )^{\frac{1-4\theta}{2\theta}}   \right ) \nonumber \\
&\leq  C P_1^2   \frac{1}{4\theta-1} t^{ \frac{4\theta -1}{2\theta }}, \quad t\geq t_0  \label{A5.24}
\end{align}
with $C = 4 \omega _1 \left ( \frac{\alpha}{2} \right )^{\frac{1-4\theta}{2\theta}} $. 
It is important to emphasize that the above estimate  holds only for $ \theta \neq \frac{1}{4} $ and we have just used it for  $ 1-4\theta <0  $ to obtain   \eqref{A5.24}.  The estimates for $A_1$ and  $A_2$ prove the desired estimate from below of lemma.
\hfill
$\Box$

As a special case one can introduce the following $\log$-order blowup result for the case of $\theta = \frac{1}{4}$.  
\begin{lem}\label{lemma4.3} Let $n=1$ and $ \theta = \frac{1}{4} $. For $ u_1 \in L^{1}({\bf R})$ the following optimal estimate holds.
$$ C_1 P_{1}^{2}\log t \leq \int _{{\bf R}} |\varphi (t,\xi)|^2 d\xi \leq C_2 P_{1}^{2} \log t , \quad t\gg 1,$$
with some constants $C_1, C_2>0$. 
\end{lem}
{\it Proof.} We consider the functions $A_1(t,\frac{1}{4})$  and $A_2(t,\frac{1}{4})$  given by \eqref{defiI_1} and \eqref{defiI_2} with $ \theta = \frac{1}{4} $.  The estimate \eqref{estimateA_1} also holds for $ \theta = \frac{1}{4} $ and it tells us the fact that 
\begin{align}\label{estimateA_1theta14}
A_1(t, \frac{1}{4}) \leq C P_1^2 , \quad t \gg 1.
\end{align}
While, by definition \eqref{defiI_2} and \eqref{A5.9} we have 
\begin{align}\label{estimateA_2theta14}
A_2(t,\frac{1}{4}) &=  \omega _1 P_1^2 \int _{\left ( \frac{\alpha}{2t} \right )^{2} }^{ \eta} e^{-\frac{2r^2}{\log(1+\sqrt{r} )}t} \left ( \frac{1 - e^{-t\frac{\log^2(1+\sqrt{r})-r^2}{\log(1+\sqrt{r})}}}{\log (1+\sqrt{r})} \right )^2 dr  \nonumber \\
&\leq 4\omega _1 P_1^2 \int _{\left ( \frac{\alpha}{2t} \right )^{2} }^{ \eta} e^{-\frac{2r^2}{\log(1+\sqrt{r})}t} \left ( \frac{1 - e^{-2t \sqrt{r}}}{\sqrt{r}} \right )^2 dr \nonumber \\
&\leq 4\omega _1 P_1^2 \int _{\left ( \frac{\alpha}{2t} \right )^{ 2}}^{ \eta} \frac{1}{r} dr \nonumber \\
&= 4\omega _1 P_1^2 \left ( \log \eta - \log \left ( \frac{\alpha}{2t} \right )^{ 2} \right )  \nonumber \\
&= 4\omega _1 P_1^2 \Big( \log \eta - 2 \log \alpha + 2\log 2 + 2 \log t  \Big) \nonumber \\
&\leq C P_1^2 \log t, \quad t\gg 1. 
\end{align}
The estimates \eqref{estimateA_1theta14} and \eqref{estimateA_2theta14} allow us to conclude  the upper estimate
\begin{align}\label{A5.27}
\int _{{\bf R}} |\varphi (t,\xi)|^2 d\xi \leq C_2 \log t, \quad t \gg 1,
\end{align}
with some constant $C_2>0$.
\vspace{0.2cm}

On the other hand, by \eqref{varphiaux} one can get
\begin{equation*}
| \varphi(t,\xi) |^2 \geq \frac{1}{2}|\varphi _1(t,\xi) |^2 - |\varphi _2(t,\xi) |^2, \quad t > 0,  \,\,\xi \in {\bf R}.
\end{equation*} 
Thus, for $ t >0$,
\begin{align}\label{lowerestimatetheta14}
\int _{{\bf R}} | \varphi(t,\xi) |^2 d\xi  &\geq \int _{t^{-1}}^{t^{-\frac{2}{3}}} | \varphi(t,\xi) |^2 d\xi \geq \frac{1}{2} \int _{t^{-1}}^{t^{-\frac{2}{3}}} |\varphi _1(t,\xi) |^2 d\xi - \int _{t^{-1}}^{t^{-\frac{2}{3}}} |\varphi _2(t,\xi) |^2 d\xi \nonumber \\
&=: P_1^2  \left ( \frac{1}{2} K_1(t) - K_2(t) \right ),
\end{align}
 where 
 \begin{align*}
 K_1(t) &:= \int _{t^{-1}}^{t^{-\frac{2}{3}}}  \frac{e^{-\frac{2 |\xi|^2}{\log(1+\sqrt{|\xi|})}t}}{\log ^2 (1+\sqrt{|\xi|})} d\xi, \\
 K_2(t) &:= \int _{t^{-1}}^{t^{-\frac{2}{3}}} \frac{e^{-2\log(1+\sqrt{|\xi|})t}}{\log^2(1+\sqrt{|\xi|})} d \xi.
 \end{align*}
We remember that
\begin{align}\label{A6.4}
\frac{1}{2}\sqrt{|\xi|} \leq  \log (1+\sqrt{|\xi|}) \leq \sqrt{|\xi|}, \quad |\xi| \leq 1. 
\end{align}
 Thus by applying Remark 3.1 one has 
\begin{align}\label{estimateK1}
 K_1(t) &= \int _{t^{-1}}^{t^{-\frac{2}{3}}}  \frac{e^{-\frac{2 |\xi|^2}{\log(1+\sqrt{|\xi|})}t}}{\log ^2 (1+\sqrt{|\xi|})} d\xi \geq \int _{t^{-1}}^{t^{-\frac{2}{3}}}  \frac{e^{-4 |\xi|^{\frac{3}{2}}  t}}{|\xi|} d\xi \nonumber \\
&= \omega _1 \int _{t^{-1}}^{t^{-\frac{2}{3}}}  \frac{e^{-4 r^{\frac{3}{2}}  t}}{r} dr  \geq  \omega _1 e^{-4} \int _{t^{-1}}^{t^{-\frac{2}{3}}}  \frac{1}{r} dr \nonumber \\
&=\omega_1 e^{-4} \big( -\frac{2}{3} \log t + \log t  \big) \nonumber \\
&= \omega_1 \frac{e^{-4}}{3} \log t , \quad t \geq 1.
\end{align}
Similarly, in the case when large $t > 1$ such that $t^{-\frac{2}{3}} < 1$ it follows from \eqref{A6.4} that
\begin{align}\label{estimateK2}
 K_2(t) &= \int _{t^{-1}}^{t^{-\frac{2}{3}}} \frac{e^{-2\log(1+\sqrt{|\xi|})t}}{\log^2(1+\sqrt{|\xi|})} d \xi \leq 4 \int _{t^{-1}}^{t^{-\frac{2}{3}}} \frac{e^{- \sqrt{\vert\xi\vert}t}}{|\xi|} d \xi \nonumber \\
&= 4 \omega _1 \int _{t^{-1}}^{t^{-\frac{2}{3}}} \frac{e^{- \sqrt{r}t}}{r} d r \leq 4 \omega _1 e^{-\sqrt{t}}\int _{t^{-1}}^{t^{-\frac{2}{3}}} \frac{1}{r} d r  \nonumber \\
&\leq \frac{4 \omega _1}{3} e^{-\sqrt{t}} \log t, \quad t \gg 1. 
\end{align}
Therefore, from \eqref{lowerestimatetheta14}, \eqref{estimateK1} and \eqref{estimateK2} one has 
\begin{align*}
\int _{{\bf R}} | \varphi(t,\xi) |^2 d\xi &\geq P_1^2  \left ( \frac{1}{2} K_1(t) - K_2(t) \right )  \\
&\geq P_1^2 \omega _1 \left (  \frac{e^{-4}}{6} \log t - \frac{4 }{3} e^{-\sqrt{t}} \log t \right )  \\
&= P_1^2 \omega _1 \log t\left (  \frac{e^{-4}}{6}  - \frac{4 }{3} e^{-\sqrt{t}} \right ) , \quad t\gg1
\end{align*}
which implies the desired estimate from below to the case $\theta=1/4$ and $n=1$:
\begin{align}\label{A5.25}
\int _{{\bf R}} | \varphi(t,\xi) |^2 d\xi \geq C P_1^2 \log t, \quad t\gg 1
\end{align}
with some constant $C > 0$. The estimates \eqref{A5.27} and \eqref{A5.25} complete the proof of lemma.
\hfill
$\Box$
\vspace{0.2cm}

Now, let us prove Theorem \ref{main-teo} at a stroke. 
\vspace{0.2cm}

{\it  Proof of Theorem \ref{main-teo}.} \,One first observes that 
\begin{equation}\label{1estimateu}
\int _{{\bf R}^n} |\hat{u}(t,\xi)|^2 d\xi \leq \int _{{\bf R}^n} |\hat{u}(t,\xi)-\varphi(t,\xi)|^2 d\xi + \int _{{\bf R}^n} |\varphi(t,\xi)|^2 d\xi.
\end{equation}
By combining Lemma \ref{lemma4.1} and Theorem \ref{maintheorem} with \eqref{1estimateu}, we have 
\begin{equation}\label{F4.47}
\int _{{\bf R}^n} | \hat{u}(t,\xi)|^2 d\xi \leq C  (P_1^2+\|u_1\|_{L^{1,2\theta}}^2) \left ( t^{-\frac{n-4\theta}{2(1-\theta)}} + \frac{1}{\theta} t^{-\frac{n-4\theta}{2\theta}}\right ) , \quad t\gg 1. 
\end{equation}

We can also observe that for $0< \theta <1/2 $ it holds that  $ 2\theta \leq 2-2\theta $. Therefore $ \frac{n-4\theta}{2-2\theta}\leq \frac{n-4\theta}{2\theta} $. Thus the  decay rate $ t^{-\frac{n-4\theta}{2\theta}} $ is faster than $ t^{-\frac{n-4\theta}{2-2\theta}} $. It results  the following upper bound to the $L^2$-norm of the Fourier transformed solution $\hat{u}$ such that
\begin{equation*}
\int _{{\bf R}^n} | \hat{u}(t,\xi)|^2 d\xi \leq C  (P_1^2+\|u_1\|_{L^{1,2\theta}}^2)(1+ \frac{1}{\theta})t^{-\frac{n-4\theta}{2(1-\theta)}} , \quad t\gg 1. 
\end{equation*}
By the Plancherel  Theorem and from \eqref{F4.47} the upper bound estimate  of  the statement of Theorem \ref{main-teo} follows with a generous constant $C>0$.  
\vspace{0.2cm} 

In order to obtain the lower bound, we observe that 
$$ |\varphi (t,\xi)| \leq |\hat{u}(t,\xi)-\varphi(t,\xi)| + |\hat{u}(t,\xi)| .$$
By Young's inequality, we may obtain 
$$ |\varphi (t,\xi)|^2 \leq 2|\hat{u}(t,\xi)-\varphi(t,\xi)|^2 + 2|\hat{u}(t,\xi)|^2 .$$
Therefore, 
\begin{align*}
 |\hat{u}(t,\xi)|^2 \geq  \frac{1}{2}|\varphi (t,\xi)|^2 - |\hat{u}(t,\xi)-\varphi(t,\xi)|^2, \quad t \geq 0 , \quad \xi \in {\bf R}^n. 
\end{align*}
Thus, 
\begin{align}\label{F4.49}
\int _{{\bf R}^n}|\hat{u}(t,\xi)|^2 d\xi \geq  \frac{1}{2}\int _{{\bf R}^n} |\varphi (t,\xi)|^2 d\xi - \int _{{\bf R}^n}|\hat{u}(t,\xi)-\varphi(t,\xi)|^2 d\xi
\end{align}

First, we consider the case $n \geq 1$ and $0< \theta \leq \frac{1}{6} $. By combining \eqref{F4.49} with the  lower estimate of Lemma \ref{lemma4.1} and estimate of Theorem \ref{maintheorem}, we obtain 
\begin{align}\label{F4.50}
\| \hat{u}(t, \cdot) \|^2 &\geq \frac{C_1}{2} P_1^2 t^{- \frac{n-4\theta}{2(1-\theta)}} -  C_2(\|u_1\|_{1}^2+\|u_1\|_{L^{1,2\theta}}^2) \left (   t^{-\frac{n}{2(1-\theta)}}  +    \frac{1}{\theta} t^{-\frac{n}{2\theta}} \right ) \nonumber \\
&= t^{- \frac{n-4\theta}{2(1-\theta)}} \left(\frac{C_1}{2} P_1^2 - 2 C_2(\|u_1\|_{1}^2+\|u_1\|_{L^{1,2\theta}}^2) \left ( t^{-\frac{4\theta}{2(1-\theta)}} + \frac{1}{\theta} t^{-\frac{n+4\theta^2 -2n\theta }{2\theta(1-\theta)}} \right )  \right),
\end{align}
for $t\gg1$ and  positive constants $C_1, C_2$. But for $ 0< \theta \leq \frac{1}{6} $ we notice that $ n+4\theta^2 -2n\theta >0 $, so that one can get
$$ \lim _{t \rightarrow \infty} \left(\frac{C_1}{2} P_1^2 - 2 C_2(\|u_1\|_{1}^2+\|u_1\|_{L^{1,2\theta}}^2) \left ( t^{-\frac{4\theta}{2(1-\theta)}} + \frac{1}{\theta} t^{-\frac{n+4\theta^2 -2n\theta }{2\theta(1-\theta)}} \right )  \right) = \frac{C_1}{2} P_1^2 .$$
Therefore, there exists $t_1>0$ such that 
\begin{align*}
\frac{C_1}{4} P_1^2 \leq \frac{C_1}{2} P_1^2 - 2 C_2(\|u_1\|_{1}^2+\|u_1\|_{L^{1,2\theta}}^2) \left ( t^{-\frac{4\theta}{2(1-\theta)}} + \frac{1}{\theta} t^{-\frac{n+4\theta^2 -2n\theta }{2\theta(1-\theta)}} \right ) \leq C_1 P_1^2, \quad t \gg t_1. 
\end{align*}
From \eqref{F4.50} it follows that 
\begin{align}\label{F4.52}
\| \hat{u}(t, \cdot) \|^2 &\geq \frac{C_1}{4} P_1^2 t^{- \frac{n-4\theta}{2(1-\theta)}} ,
\end{align}
for $ t\gg1 $. 

The estimate \eqref{F4.52} implies the desired estimate for lower bound in $t$ in the case when $n \geq 1$ and $0< \theta \leq \frac{1}{6} $.

Analogously, we may obtain the results for $ n=1 $ with $ \frac{1}{6}< \theta< \frac{1}{4} $, and for $n \geq 2$ with $ \frac{1}{6} <\theta \leq \frac{5}{12} $ based on the results of  Theorem \ref{maintheorem} for these values of $\theta$ and Lemma \ref{lemma4.1}.  

This completes the proof of Theorem \ref{main-teo}.
\hfill
$\Box$
\vspace{0.4cm}


{\it{Proof of Theorem \ref{main-teo-2}.}} The proof of Theorem \ref{main-teo-2} can be obtained in the same way as in Theorem \ref{main-teo}, but using Lemmas \ref{lemma4.2} and \ref{lemma4.3} instead of Lemma \ref{lemma4.1} and observing that the estimates to $ \|\varphi(t, \cdot)\|^2 $ in Lemmas \ref{lemma4.2} and \ref{lemma4.3} are also worse than the estimates to $ \|\hat{u}(t,\cdot) - \varphi(t,\cdot)\|^2 $ in Theorem \ref{maintheorem}. 
\hfill
$\Box$

\par
\vspace{0.5cm}
\noindent{\em Acknowledgements.}
The work of the first author (A. PISKE) was financed in part by the Coordena\c{c}\~ao de Aperfei\c{c}oamento de Pessoal de N\'ivel Superior - Brasil (Capes) -  Finance Code 001, and the work of the second author (R. C. CHAR\~AO) was partially supported by PRINT/CAPES- Process 88881.310536/2018-00. The work of the third author (R. IKEHATA) was supported in part by Grant-in-Aid for Scientific Research (C) 20K03682  of JSPS.




\end{document}